\documentclass[numbers=enddot,12pt,final,onecolumn,notitlepage]{scrartcl}%
\usepackage[headsepline,footsepline,manualmark]{scrlayer-scrpage}
\usepackage[all,cmtip]{xy}
\usepackage{amssymb}
\usepackage{amsmath}
\usepackage{amsthm}
\usepackage{framed}
\usepackage{comment}
\usepackage{color}
\usepackage[breaklinks=true]{hyperref}
\usepackage[sc]{mathpazo}
\usepackage[T1]{fontenc}
\usepackage{needspace}
\usepackage{tabls}
\usepackage{tikz}
\providecommand{\U}[1]{\protect\rule{.1in}{.1in}}
\usetikzlibrary{arrows.meta}
\theoremstyle{definition}
\newtheorem{theo}{Theorem}[section]
\newenvironment{theorem}[1][]
{\begin{theo}[#1]\begin{leftbar}}
{\end{leftbar}\end{theo}}
\newtheorem{lem}[theo]{Lemma}
\newenvironment{lemma}[1][]
{\begin{lem}[#1]\begin{leftbar}}
{\end{leftbar}\end{lem}}
\newtheorem{prop}[theo]{Proposition}
\newenvironment{proposition}[1][]
{\begin{prop}[#1]\begin{leftbar}}
{\end{leftbar}\end{prop}}
\newtheorem{defi}[theo]{Definition}
\newenvironment{definition}[1][]
{\begin{defi}[#1]\begin{leftbar}}
{\end{leftbar}\end{defi}}
\newtheorem{remk}[theo]{Remark}
\newenvironment{remark}[1][]
{\begin{remk}[#1]\begin{leftbar}}
{\end{leftbar}\end{remk}}
\newtheorem{coro}[theo]{Corollary}
\newenvironment{corollary}[1][]
{\begin{coro}[#1]\begin{leftbar}}
{\end{leftbar}\end{coro}}
\newtheorem{conv}[theo]{Convention}

\newtheorem{quest}[theo]{Question}
\newenvironment{question}[1][]
{\begin{quest}[#1]\begin{leftbar}}
{\end{leftbar}\end{quest}}
\newtheorem{warn}[theo]{Warning}

\newtheorem{conj}[theo]{Conjecture}

\newtheorem{exam}[theo]{Example}
\newenvironment{example}[1][]
{\begin{exam}[#1]\begin{leftbar}}
{\end{leftbar}\end{exam}}

\newcommand{\arinj}{\ar@{_{(}->}}
\newcommand{\arinjrev}{\ar@{^{(}->}}
\newcommand{\arsurj}{\ar@{->>}}
\newcommand{\arelem}{\ar@{|->}}
\newcommand{\arback}{\ar@{<-}}
\DeclareMathOperator{\sep}{sep}
\DeclareMathOperator{\Sep}{Sep}

\let\sumnonlimits\sum
\let\prodnonlimits\prod
\let\cupnonlimits\bigcup
\let\capnonlimits\bigcap
\let\bigoplusnonlimits\bigoplus
\renewcommand{\sum}{\sumnonlimits\limits}
\renewcommand{\prod}{\prodnonlimits\limits}
\renewcommand{\bigcup}{\cupnonlimits\limits}
\renewcommand{\bigcap}{\capnonlimits\limits}
\renewcommand{\bigoplus}{\bigoplusnonlimits\limits}

\excludecomment{noncompile}
\ihead{The $q$-Poincar\'e sum}\ohead{page \thepage}
\begin{document}

\title{The $q$-Poincar\'{e} sum and separability in a Coxeter group}
\author{Ben Adenbaum, Darij Grinberg}
\date{November 4, 2026}
\maketitle

\begin{abstract}
\textbf{Abstract.} Given a finite Coxeter group $W$ and an element $q$ of a
commutative ring $\mathbf{k}$, we define the element $\mathbf{L}_{q}%
:=\sum_{w\in W}q^{\ell\left(  w\right)  }w$ of its group algebra
$\mathbf{k}\left[  W\right]  $. We say that an element $w\in W$ is
$W$\emph{-length-symmetric} if its coefficient in the commutator $\left[
\mathbf{L}_{a},\mathbf{L}_{b}\right]  $ is $0$ for all $\mathbf{k}$ and all
$a,b\in\mathbf{k}$. We prove several sufficient criteria for elements to be
$W$-length-symmetric. In particular, we show that any separable permutation
(i.e., permutation avoiding the patterns $2413$ and $3142$) in the symmetric
group $S_{n}$ is $S_{n}$-length-symmetric. More generally, we construct a
notion of separable elements of a finite Coxeter group $W$ (which, for finite
Weyl groups, agrees with the eponymous notion introduced by Gaetz and Gao),
and show that all such elements are $W$-length-symmetric. This is only a
sufficient condition; other $W$-length-symmetric elements include (twisted)
involutions and elements of dihedral parabolic subgroups.

\end{abstract}

Fix a commutative ring $\mathbf{k}$. Consider the $n$-th symmetric group
$S_{n}$ for a given nonnegative integer $n$. Any permutation $w\in S_{n}$ has
a \emph{length} $\ell\left(  w\right)  $, which is defined as the number of
inversions of $w$ (that is, of pairs $\left(  i,j\right)  $ of integers with
$1\leq i<j\leq n$ and $w\left(  i\right)  >w\left(  j\right)  $). Using these
lengths, we can define an element%
\[
\mathbf{L}_{q}:=\sum_{w\in S_{n}}q^{\ell\left(  w\right)  }w\in\mathbf{k}%
\left[  S_{n}\right]
\]
in the group algebra $\mathbf{k}\left[  S_{n}\right]  $ of $S_{n}$, for any
fixed scalar $q\in\mathbf{k}$. For instance, $\mathbf{L}_{1}$ and
$\mathbf{L}_{-1}$ are two central elements of $\mathbf{k}\left[  S_{n}\right]
$. (When $n!$ is invertible in $\mathbf{k}$, they equal $n!$ times the central
idempotents corresponding to the trivial and sign representations of $S_{n}$,
respectively.) For general $q$, however, the element $\mathbf{L}_{q}$ is not
central, and does not even necessarily commute with another such element
$\mathbf{L}_{r}$ for $r\neq q$.

However, it turns out that if we fix two scalars $q,r\in\mathbf{k}$, then many
permutations $w\in S_{n}$ do not appear in the commutator $\left[
\mathbf{L}_{q},\mathbf{L}_{r}\right]  :=\mathbf{L}_{q}\mathbf{L}%
_{r}-\mathbf{L}_{r}\mathbf{L}_{q}$. Namely, if $w\in S_{n}$ is a
\emph{separable} permutation, then the coefficient of $w$ in $\left[
\mathbf{L}_{q},\mathbf{L}_{r}\right]  $ is $0$. The notion of a separable
permutation can be defined in many ways, e.g., algorithmically or by the
avoidance of two patterns ($2413$ and $3142$); we refer to Ghys's book
\cite[\textquotedblleft Separable permutations\textquotedblright]{Ghys16} for
four equivalent definitions and a survey of their rich theory.

In this paper, we will prove the fact just mentioned, and essentially
generalize it to arbitrary finite Coxeter groups.

\section{The theorems}

We assume the basic properties of Coxeter groups (see \cite{BjoBre},
\cite[Chapters 3--4]{Davis}, \cite[\S 1, \S 2, \S 9]{Lusztig}, \cite[\S IV.1]%
{Bourbaki-Lie4-6}, \cite[\S 4]{Heckman}, \cite{Digne}, \cite{Hiller},
\cite[Chapter 5]{Humphreys}). A Coxeter system $\left(  W,S\right)  $ consists
of a Coxeter group $W$ and its set $S$ of simple reflections. If $\left(
W,S\right)  $ is a Coxeter system, and $w\in W$ is any element, then
$\ell\left(  w\right)  $ denotes the Coxeter length of $w$ (that is, the
smallest length of a list of simple reflections having product $w$). Note that
$\ell\left(  w^{-1}\right)  =\ell\left(  w\right)  $ for each $w\in W$.

If $W$ is a finite Coxeter group, then there exists a unique element $w\in W$
whose length is maximum (see \cite[\S 11.1]{Lusztig}); we shall call it
$w_{0}$. This element $w_{0}$ is always an involution (i.e., satisfies
$w_{0}^{2}=\operatorname*{id}$), since its inverse $w_{0}^{-1}$ must also have
maximum length. It is known as the \emph{longest element} or the \emph{longest
word} of $W$.

\subsection{$W$-length-symmetry}

We now define our main object of study:

\begin{definition}
\label{def.Lq}Let $\left(  W,S\right)  $ be a Coxeter system with finite
Coxeter group $W$. Let $\mathbf{k}$ be a commutative ring. For any
$q\in\mathbf{k}$, we define the element%
\[
\mathbf{L}_{q}\left(  W\right)  :=\sum_{w\in W}q^{\ell\left(  w\right)  }w
\]
in the group algebra $\mathbf{k}\left[  W\right]  $. We denote this element
simply by $\mathbf{L}_{q}$ if the group $W$ is clear from the context. We call
it the $q$\emph{-Poincar\'{e} sum} of $W$.
\end{definition}

\begin{example}
\label{exa.LqSn}Let $n$ be a positive integer. The symmetric group $S_{n}$ is
the Coxeter group of type $A_{n-1}$; the corresponding Coxeter system is
$\left(  S_{n},S\right)  $, where $S=\left\{  s_{1},s_{2},\ldots
,s_{n-1}\right\}  $ is the set of all simple reflections in $S_{n}$ (so that
$s_{i}$ is the transposition swapping $i$ with $i+1$). The length $\ell\left(
w\right)  $ of a permutation $w\in S_{n}$ is well-known to equal the number of
inversions of $w$ (that is, of pairs $\left(  i,j\right)  $ of elements of
$\left\{  1,2,\ldots,n\right\}  $ satisfying $i<j$ and $w\left(  i\right)
>w\left(  j\right)  $). The $q$-Poincar\'{e} sum $\mathbf{L}_{q}\left(
S_{n}\right)  =\sum_{w\in S_{n}}q^{\ell\left(  w\right)  }w\in\mathbf{k}%
\left[  S_{n}\right]  $ is then a linear combination of all permutations of
$\left\{  1,2,\ldots,n\right\}  $, where each $w$ appears with coefficient
$q^{\ell\left(  w\right)  }$. For instance,%
\[
\mathbf{L}_{q}\left(  S_{3}\right)  =\operatorname*{id}+\,qs_{1}+qs_{2}%
+q^{2}s_{1}s_{2}+q^{2}s_{2}s_{1}+q^{3}w_{0},
\]
where $w_{0}=s_{1}s_{2}s_{1}=s_{2}s_{1}s_{2}$ is the permutation sending
$1,2,3$ to $3,2,1$.

If $q$ is a nonnegative real, then the element $\mathbf{L}_{q}\left(
S_{n}\right)  $ can be viewed as a measure on $S_{n}$ (by treating each
coefficient $q^{\ell\left(  w\right)  }$, divided by the total sum of these
coefficients, as the measure of $w\in S_{n}$). This measure is known as the
\emph{Mallows measure} $\operatorname*{Mallows}\left(  n,q\right)  $ (see,
e.g., \cite{Korotk16}).
\end{example}

\begin{example}
Let $\left(  W,S\right)  $ be any Coxeter system with finite Coxeter group
$W$. Then,
\[
\mathbf{L}_{1}\left(  W\right)  =\sum_{w\in W}w\ \ \ \ \ \ \ \ \ \ \text{and}%
\ \ \ \ \ \ \ \ \ \ \mathbf{L}_{-1}\left(  W\right)  =\sum_{w\in
W}\underbrace{\left(  -1\right)  ^{\ell\left(  w\right)  }}%
_{=\operatorname*{sign}w}w
\]
are two elements in the center of $\mathbf{k}\left[  W\right]  $, and in fact
(up to a factor of $\left\vert W\right\vert $) are the central idempotents
corresponding to the trivial and sign representations of $W$, respectively
(when $\left\vert W\right\vert $ is invertible). Of course, $\mathbf{L}%
_{0}\left(  W\right)  =1$ lies in the center of $\mathbf{k}\left[  W\right]
$, too.
\end{example}

However, in general, $\mathbf{L}_{q}=\mathbf{L}_{q}\left(  W\right)  $ does
not lie in the center of $\mathbf{k}\left[  W\right]  $ when $q\notin\left\{
1,-1,0\right\}  $; and in fact, the elements $\mathbf{L}_{q}$ for two
different values of $q$ usually do not commute. That is, most finite Coxeter
groups $W$ do \textbf{not} satisfy the equality $\mathbf{L}_{a}\mathbf{L}%
_{b}=\mathbf{L}_{b}\mathbf{L}_{a}$ when $a$ and $b$ are generic.

This note is devoted to salvaging this equality to some extent. First, we
shall show that it holds when $W$ is dihedral (or smaller):

\begin{theorem}
\label{thm.comm.dihe}Let $\left(  W,S\right)  $ be a Coxeter system, where $W$
is finite and $\left\vert S\right\vert \leq2$. Then,
\[
\mathbf{L}_{a}\mathbf{L}_{b}=\mathbf{L}_{b}\mathbf{L}_{a}
\]
for any $a,b\in\mathbf{k}$.
\end{theorem}

\begin{example}
Consider the symmetric group $S_{3}$ as the Coxeter group $A_{2}$. Then one
can observe that
\begin{align*}
&  \mathbf{L}_{a}\mathbf{L}_{b}\\
&  =(1+2ab+2a^{2}b^{2}+a^{3}b^{3})\operatorname*{id}+(a+b+ab^{2}+a^{2}%
b+a^{3}b^{2}+a^{2}b^{3})(s_{1}+s_{2})\\
&  \ \ \ \ +(ab+a^{2}+b^{2}+a^{2}b^{2}+ab^{3}+a^{3}b)(s_{1}s_{2}+s_{2}%
s_{1})+(2ab^{2}+2a^{2}b+a^{3}+b^{3})w_{0}%
\end{align*}
with the notations of Example \ref{exa.LqSn}. This is visibly symmetric in $a$
and $b$, so it confirms Theorem \ref{thm.comm.dihe}.
\end{example}

We will prove Theorem \ref{thm.comm.dihe}, and many other results we shall
soon state, in the next section (Section \ref{sec.pfs}).

Next, even though the products $\mathbf{L}_{a}\mathbf{L}_{b}$ and
$\mathbf{L}_{b}\mathbf{L}_{a}$ are usually distinct for general $W$, we shall
show that certain coefficients of these products are equal. For this, we
introduce a notation:

\begin{definition}
\label{def.coeff}Let $W$ be a group. Let $a\in\mathbf{k}\left[  W\right]  $ be
an element of its group algebra. Let $w\in W$. Then, the $w$-coefficient of
$a$ will be denoted by $\left[  w\right]  a$.
\end{definition}

\begin{definition}
\label{def.Wls}Let $\left(  W,S\right)  $ be a Coxeter system with finite
Coxeter group $W$. Let $w\in W$. We say that the element $w$ is $W$%
\emph{-length-symmetric} if every commutative ring $\mathbf{k}$ and every
$a,b\in\mathbf{k}$ satisfy%
\begin{equation}
\left[  w\right]  \left(  \mathbf{L}_{a}\mathbf{L}_{b}\right)  =\left[
w\right]  \left(  \mathbf{L}_{b}\mathbf{L}_{a}\right)  . \label{eq.def.Wls.eq}%
\end{equation}

\end{definition}

Alternatively, $W$-length-symmetry can be characterized combinatorially as follows:

\begin{proposition}
\label{prop.Wls.comb}Let $\left(  W,S\right)  $ be a Coxeter system with
finite Coxeter group $W$. Let $w\in W$. Then, $w$ is $W$-length-symmetric if
and only if for every $p,q\in\mathbb{N}$, we have%
\begin{align*}
&  \left(  \text{number of }u\in W\text{ satisfying }\ell\left(  u\right)
=p\text{ and }\ell\left(  u^{-1}w\right)  =q\right) \\
&  =\left(  \text{number of }u\in W\text{ satisfying }\ell\left(  u\right)
=p\text{ and }\ell\left(  wu^{-1}\right)  =q\right)  .
\end{align*}

\end{proposition}

\begin{proof}
By Definition \ref{def.Lq}, we have $\mathbf{L}_{a}=\sum_{w\in W}%
a^{\ell\left(  w\right)  }w$ and $\mathbf{L}_{b}=\sum_{w\in W}b^{\ell\left(
w\right)  }w$ for all $a,b\in\mathbf{k}$. Thus, for all $a,b\in\mathbf{k}$, we
have
\begin{align*}
\mathbf{L}_{a}\mathbf{L}_{b}  &  =\left(  \sum_{w\in W}a^{\ell\left(
w\right)  }w\right)  \left(  \sum_{w\in W}b^{\ell\left(  w\right)  }w\right)
=\sum_{\left(  u,v\right)  \in W\times W}a^{\ell\left(  u\right)  }%
b^{\ell\left(  v\right)  }uv\ \ \ \ \ \ \ \ \ \ \text{and}\\
\mathbf{L}_{b}\mathbf{L}_{a}  &  =\left(  \sum_{w\in W}b^{\ell\left(
w\right)  }w\right)  \left(  \sum_{w\in W}a^{\ell\left(  w\right)  }w\right)
=\sum_{\left(  v,u\right)  \in W\times W}b^{\ell\left(  v\right)  }%
a^{\ell\left(  u\right)  }vu\\
&  =\sum_{\left(  v,u\right)  \in W\times W}a^{\ell\left(  u\right)  }%
b^{\ell\left(  v\right)  }vu.
\end{align*}
Hence, each $w\in W$ and $a,b\in\mathbf{k}$ satisfy%
\[
\left[  w\right]  \left(  \mathbf{L}_{a}\mathbf{L}_{b}\right)  =\sum
_{\substack{\left(  u,v\right)  \in W\times W;\\uv=w}}a^{\ell\left(  u\right)
}b^{\ell\left(  v\right)  }=\sum_{u\in W}a^{\ell\left(  u\right)  }%
b^{\ell\left(  u^{-1}w\right)  }%
\]
(here, we have substituted $\left(  u,u^{-1}w\right)  $ for $\left(
u,v\right)  $) and%
\[
\left[  w\right]  \left(  \mathbf{L}_{b}\mathbf{L}_{a}\right)  =\sum
_{\substack{\left(  v,u\right)  \in W\times W;\\vu=w}}a^{\ell\left(  u\right)
}b^{\ell\left(  v\right)  }=\sum_{u\in W}a^{\ell\left(  u\right)  }%
b^{\ell\left(  wu^{-1}\right)  }%
\]
(here, we have substituted $\left(  wu^{-1},u\right)  $ for $\left(
v,u\right)  $). Hence, the equality (\ref{eq.def.Wls.eq}) can be rewritten as%
\begin{equation}
\sum_{u\in W}a^{\ell\left(  u\right)  }b^{\ell\left(  u^{-1}w\right)  }%
=\sum_{u\in W}a^{\ell\left(  u\right)  }b^{\ell\left(  wu^{-1}\right)  }.
\label{eq.def.Wls.eq2}%
\end{equation}
This equality is a polynomial identity in $a$ and $b$. Thus, it holds for all
commutative rings $\mathbf{k}$ and all $a,b\in\mathbf{k}$ if and only if it
holds as a polynomial identity, i.e., if and only if we have%
\begin{equation}
\sum_{u\in W}X^{\ell\left(  u\right)  }Y^{\ell\left(  u^{-1}w\right)  }%
=\sum_{u\in W}X^{\ell\left(  u\right)  }Y^{\ell\left(  wu^{-1}\right)  }
\label{eq.def.Wls.eq3}%
\end{equation}
in the polynomial ring $\mathbb{Z}\left[  X,Y\right]  $ (because
(\ref{eq.def.Wls.eq2}) can be derived from (\ref{eq.def.Wls.eq3}) by
evaluating both sides at $X=a$ and $Y=b$, whereas (\ref{eq.def.Wls.eq3}) is
the particular case of (\ref{eq.def.Wls.eq2}) for $\mathbf{k}=\mathbb{Z}%
\left[  X,Y\right]  $ and $a=X$ and $b=Y$). But the identity
(\ref{eq.def.Wls.eq3}), in turn, is equivalent to saying that for every
$p,q\in\mathbb{N}$, we have%
\begin{align}
&  \left(  \text{number of }u\in W\text{ satisfying }\ell\left(  u\right)
=p\text{ and }\ell\left(  u^{-1}w\right)  =q\right) \nonumber\\
&  =\left(  \text{number of }u\in W\text{ satisfying }\ell\left(  u\right)
=p\text{ and }\ell\left(  wu^{-1}\right)  =q\right)  \label{eq.def.Wls.eq4}%
\end{align}
(since the former number is the $X^{p}Y^{q}$-coefficient on the left-hand side
of (\ref{eq.def.Wls.eq3}), whereas the latter number is the $X^{p}Y^{q}%
$-coefficient on the right-hand side of (\ref{eq.def.Wls.eq3})).

Combining these equivalences, we conclude that the equality
(\ref{eq.def.Wls.eq}) holds for all commutative rings $\mathbf{k}$ and all
$a,b\in\mathbf{k}$ if and only if the equality (\ref{eq.def.Wls.eq4}) holds
for every $p,q\in\mathbb{N}$. But the former equality is exactly the
definition of \textquotedblleft$W$-length-symmetric\textquotedblright; thus,
we have shown that $w$ is $W$-length-symmetric if and only if the equality
(\ref{eq.def.Wls.eq4}) holds for every $p,q\in\mathbb{N}$. This proves
Proposition \ref{prop.Wls.comb}.
\end{proof}

\begin{remark}
The characterization of $W$-length-symmetry in Proposition \ref{prop.Wls.comb}
makes perfect sense even if $W$ is infinite, as long as $S$ is finite (since
the number of $u\in W$ satisfying $\ell\left(  u\right)  =p$ is still finite
in this case). However, we will not concern ourselves with infinite Coxeter
groups here.
\end{remark}

\begin{remark}
\label{rmk.Wls.generic}As we have already noticed, the equality
(\ref{eq.def.Wls.eq2}) is a polynomial identity in $a$ and $b$. In other
words, the equality (\ref{eq.def.Wls.eq}) is a polynomial identity in $a$ and
$b$. Thus, in order for it to hold for all $\mathbf{k}$ and $a,b$, it suffices
that it hold when $\mathbf{k}=\mathbb{Q}$ and $a$ and $b$ are positive
integers (since two polynomials over $\mathbb{Q}$ are identical if they agree
on all positive integers). Hence, an element $w\in W$ is $W$-length-symmetric
if and only if it satisfies (\ref{eq.def.Wls.eq}) for $\mathbf{k}=\mathbb{Q}$
and for all positive integers $a$ and $b$. This gives yet another
combinatorial way to think of $W$-length-symmetry.
\end{remark}

\begin{example}
Consider the symmetric group $S_{n}$ as a Coxeter group. It is easy to see
(and follows from Corollary \ref{cor.Wls.dih} below) that each element of
$S_{n}$ is $S_{n}$-length-symmetric for all $n\leq3$. For $n=4$, all elements
of $S_{4}$ are $S_{4}$-length-symmetric except for the two permutations with
one-line notations $\left(  2,4,1,3\right)  $ and $\left(  3,1,4,2\right)  $.
These two permutations are precisely the two patterns that the \emph{separable
permutations} are defined to avoid. And we shall indeed see (Theorem
\ref{thm.separable_equivalence} and Theorem \ref{thm.sep}) that every
separable permutation $w\in S_{n}$ is $S_{n}$-length-symmetric, although the
converse is not true. (We shall also see that every involution $w\in S_{n}$ is
$S_{n}$-length-symmetric, but not every involution is separable.)
\end{example}

\begin{remark}
\label{rmk.countSn}One might wonder how many $W$-length-symmetric elements
exist for a Coxeter group $W$. For instance, for each $n\in\mathbb{N}$, let
$\omega_{n}$ be the number of $S_{n}$-length-symmetric elements of the Coxeter
group $S_{n}$. Here are the first few values:%
\[%
\begin{tabular}
[c]{|c||c|c|c|c|c|c|c|c|c|}\hline
$n$ & $0$ & $1$ & $2$ & $3$ & $4$ & $5$ & $6$ & $7$ & $8$\\\hline
$\omega_{n}$ & $1$ & $1$ & $2$ & $6$ & $22$ & $94$ & $424$ & $2032$ &
$9980$\\\hline
\end{tabular}
\ \ \ \
\]
This sequence is not in the OEIS as of 2026-11-03, even without $\omega_{0}$.

Somewhat related is \href{https://oeis.org/A180389}{OEIS sequence A180389},
which counts the permutations $w\in S_{n}$ satisfying $\operatorname*{des}%
w=\operatorname*{des}\left(  w^{-1}\right)  $, where $\operatorname*{des}u$
denotes the number of descents of $u$ (that is, of simple reflections $s\in S$
satisfying $\ell\left(  us\right)  <\ell\left(  u\right)  $). Indeed, any
$W$-length-symmetric element $w\in S_{n}$ must satisfy $\operatorname*{des}%
w=\operatorname*{des}\left(  w^{-1}\right)  $, as can easily be seen from the
$p=1$ part of Proposition \ref{prop.Wls.comb}; but the converse is quite far
from truth.
\end{remark}

We shall now give several criteria for when elements of a Coxeter group are
$W$-length-symmetric. The first one is an easy consequence of Theorem
\ref{thm.comm.dihe}:

\begin{corollary}
\label{cor.Wls.dih}Let $\left(  W,S\right)  $ be a Coxeter system, where $W$
is finite and $\left\vert S\right\vert \leq2$. Then, each $w\in W$ is $W$-length-symmetric.
\end{corollary}

We shall now find $W$-length-symmetric elements in arbitrary finite Coxeter
groups. Some are fairly simple:

\begin{proposition}
\label{prop.Wls.inv}Let $\left(  W,S\right)  $ be a Coxeter system, where $W$
is finite. Then, each involution $w\in W$ (that is, each $w\in W$ satisfying
$w^{2}=1$) is $W$-length-symmetric.
\end{proposition}

\begin{proposition}
\label{prop.Wls.len2}Let $\left(  W,S\right)  $ be a Coxeter system, where $W$
is finite. Let $w\in W$ satisfy $\ell\left(  w\right)  \leq2$. Then, $w$ is
$W$-length-symmetric.
\end{proposition}

\begin{proposition}
\label{prop.Wls.inverse}Let $\left(  W,S\right)  $ be a Coxeter system, where
$W$ is finite. Let $w\in W$. Then, $w$ is $W$-length-symmetric if and only if
$w^{-1}$ is $W$-length-symmetric.
\end{proposition}

\begin{proposition}
\label{prop.Wls.ref}Let $\left(  W,S\right)  $ be a Coxeter system, where $W$
is finite. Let $w\in W$. Then, $w$ is $W$-length-symmetric if and only if
$w_{0}w$ is $W$-length-symmetric.
\end{proposition}

We can generalize Proposition \ref{prop.Wls.inv} from involutions to
\textquotedblleft twisted involutions\textquotedblright:

\begin{proposition}
\label{prop.Wls.twinv}Let $\left(  W,S\right)  $ be a Coxeter system, where
$W$ is finite. Let $\phi:W\rightarrow W$ be an automorphism of Coxeter groups
(that is, $\phi\left(  S\right)  =S$). Let $w\in W$ be such that $w^{-1}%
=\phi\left(  w\right)  $. Then, $w$ is $W$-length-symmetric.
\end{proposition}

\begin{example}
Let $\left(  W,S\right)  $ be a Coxeter system, where $W$ is finite. Then,
conjugation by $w_{0}$ (that is, the map $W\rightarrow W,\ w\mapsto
w_{0}ww_{0}^{-1}$) is an automorphism of Coxeter groups. The elements $w\in W$
satisfying $w^{-1}=w_{0}ww_{0}^{-1}$ are sometimes called the \emph{twisted
involutions}. Proposition \ref{prop.Wls.twinv} shows that they are
$W$-length-symmetric. (That said, this also follows from Proposition
\ref{prop.Wls.inv} and Proposition \ref{prop.Wls.ref}, since $w$ is a twisted
involution if and only if $w_{0}w$ is an involution.)
\end{example}

Next, we shall show that $W$-length-symmetry does not really depend on $W$, in
the sense that an element of a standard parabolic subgroup $W_{I}$ of $W$ is
$W$-length-symmetric if and only if it is $W_{I}$-length-symmetric. First we
recall the definition of standard parabolic subgroups:

\begin{definition}
Let $\left(  W,S\right)  $ be a Coxeter system. If $I$ is a subset of $S$,
then $W_{I}$ denotes the subgroup of $W$ generated by $I$. This is called a
\emph{standard parabolic subgroup} of $W$. It is known (see, e.g.,
\cite[Proposition 9.5]{Lusztig}) that $W_{I}$ is again a Coxeter group, part
of the Coxeter system $\left(  W_{I},I\right)  $.
\end{definition}

\begin{proposition}
\label{prop.Wls.WI}Let $\left(  W,S\right)  $ be a Coxeter system, where $W$
is finite. Let $I$ be a subset of $S$. Let $w\in W_{I}$. Then, $w$ is
$W$-length-symmetric if and only if $w$ is $W_{I}$-length-symmetric.
\end{proposition}

Next, we discuss $W$-length-symmetry for reducible Coxeter groups. A Coxeter
group is said to be \emph{reducible} if it has an orthogonal decomposition;
this is defined as follows:

\begin{definition}
\label{def.orthdec} Let $\left(  W,S\right)  $ be a Coxeter system. Let
$m\left(  i,j\right)  $ denote the order of $ij$ in $W$ for all $i,j\in S$. An
\emph{orthogonal decomposition} of $\left(  W,S\right)  $ means a pair
$\left(  I,J\right)  $ of two disjoint nonempty subsets $I$ and $J$ of $S$
such that $S=I\cup J$ and such that all $i\in I$ and $j\in J$ satisfy
$m\left(  i,j\right)  =2$ (thus $ij=ji$). We shall denote this decomposition
as $S=I\oplus J$.
\end{definition}

It is easy to see that if $S=I\oplus J$ is an orthogonal decomposition, then
$W=W_{I}W_{J}\cong W_{I}\times W_{J}$.

\begin{proposition}
\label{prop.Wls.I+J}Let $\left(  W,S\right)  $ be a Coxeter system, where $W$
is finite. Let $S=I\oplus J$ be an orthogonal decomposition. Let $u\in W_{I}$
and $v\in W_{J}$. Assume that $u$ is $W_{I}$-length-symmetric and $v$ is
$W_{J}$-length-symmetric. Then, $uv$ is $W$-length-symmetric.
\end{proposition}

\begin{remark}
The converse of Proposition \ref{prop.Wls.I+J} is false. A large source of
counterexamples comes from Proposition \ref{prop.Wls.twinv}. Namely, let
$\left(  W,S\right)  $ be a Coxeter system that is a \textquotedblleft direct
product\textquotedblright\ of two isomorphic Coxeter systems -- i.e., that has
an orthogonal decomposition $S=I\oplus J$ and an isomorphism $\rho:\left(
W_{I},I\right)  \rightarrow\left(  W_{J},J\right)  $ of Coxeter groups. Let
$\phi:\left(  W,S\right)  \rightarrow\left(  W,S\right)  $ be the automorphism
of Coxeter groups that is induced by $\rho$ (that is, sends $I$ to $J$ via
$\rho$, and sends $J$ to $I$ via $\rho^{-1}$). Fix any $u\in W_{I}$, and let
$v:=\rho\left(  u^{-1}\right)  \in W_{J}$. Then, the element $uv\in W$
satisfies $\left(  uv\right)  ^{-1}=\phi\left(  uv\right)  $ (since the
definition of $\phi$ yields $\phi\left(  uv\right)  =\underbrace{\rho\left(
u\right)  }_{=v^{-1}}\underbrace{\rho^{-1}\left(  v\right)  }_{=u^{-1}}%
=v^{-1}u^{-1}=\left(  uv\right)  ^{-1}$), and thus is $W$-length-symmetric by
Proposition \ref{prop.Wls.twinv}. But $u$ is not necessarily $W_{I}%
$-length-symmetric, since it can be any element of $W_{I}$.
\end{remark}

\begin{remark}
\label{rmk.not-enough}The above propositions provide several sufficient
criteria for $W$-length-symmetry; but even in combination, these criteria do
not cover all $W$-length-symmetric elements. For example, the permutation
$w\in S_{7}$ with one-line notation $\left(  2,3,5,1,7,6,4\right)  $ is
$S_{7}$-length-symmetric (checked by direct computation), but is neither an
involution nor a twisted involution, nor belongs to a subgroup $W_{I}$ for a
proper subset $I$ of $S$, nor becomes an element of such a subgroup $W_{I}$
after multiplication by $w_{0}$.
\end{remark}

\begin{remark}
\label{rmk.hecke}One can also define a variant of $W$-length-symmetry in the
Hecke algebra. Let $\left(  W,S\right)  $ be a Coxeter system with $W$ finite,
and let $\mathcal{H}$ be its Hecke algebra (see, e.g., \cite[\S 7.4]%
{Humphreys} or \cite{Digne}) over the commutative ring $\mathbf{k}$ with a
parameter $q\in\mathbf{k}$. For any $r\in\mathbf{k}$, define
\[
\mathbf{L}_{r}^{\mathcal{H}}:=\sum_{w\in W}r^{\ell\left(  w\right)  }T_{w}%
\in\mathcal{H}.
\]
Then, an element $w\in W$ can be called \emph{Hecke-}$W$%
\emph{-length-symmetric} if the coefficient of $T_{w}$ in $\mathbf{L}%
_{a}^{\mathcal{H}}\mathbf{L}_{b}^{\mathcal{H}}$ equals the coefficient of
$T_{w}$ in $\mathbf{L}_{b}^{\mathcal{H}}\mathbf{L}_{a}^{\mathcal{H}}$ for
arbitrary $\mathbf{k}$ and $q,a,b\in\mathbf{k}$. This is a stronger property
than $W$-length-symmetry; in particular, the permutations $\left[
2431\right]  ,\ \left[  3241\right]  ,\ \left[  4132\right]  ,\ \left[
4213\right]  $ in $S_{4}$ (written in one-line notation) are $S_{4}%
$-length-symmetric but not Hecke-$S_{4}$-length-symmetric.

Finding Hecke analogues of the above results on $W$-length-symmetry appears
doable but not always trivial. An easy first step is generalizing Theorem
\ref{thm.comm.dihe}: If $\left\vert S\right\vert \leq2$, then%
\begin{equation}
\mathbf{L}_{a}^{\mathcal{H}}\mathbf{L}_{b}^{\mathcal{H}}=\mathbf{L}%
_{b}^{\mathcal{H}}\mathbf{L}_{a}^{\mathcal{H}}\ \ \ \ \ \ \ \ \ \ \text{for
all }q,a,b\in\mathbf{k}. \label{eq.thm.comm.dihe.hecke}%
\end{equation}
This can be proved in roughly the same way as we prove Theorem
\ref{thm.comm.dihe} below (viz., by replacing each $w$ by $T_{w}$ in Lemma
\ref{lem.dih.sti}, and slightly modifying its proof). Furthermore, Proposition
\ref{prop.Wls.inv}, Proposition \ref{prop.Wls.inverse} and Proposition
\ref{prop.Wls.twinv} also remain true if we replace \textquotedblleft%
$W$-length-symmetric\textquotedblright\ by \textquotedblleft Hecke-$W$%
-length-symmetric\textquotedblright, and the proof we give below generalizes
to this as well (since $\mathcal{H}$ also has an \textquotedblleft
antipode\textquotedblright, i.e., an involutive $\mathbf{k}$-algebra
anti-automorphism that sends each $T_{w}$ to $T_{w^{-1}}$). On the other hand,
a version of Proposition \ref{prop.Wls.ref} appears to be missing for
$\mathcal{H}$. The generalizability of Proposition \ref{prop.Wls.WI} is an
interesting question.
\end{remark}

\subsection{$W$-separability}

We now turn to generalizing separable permutations. As we recall, any finite
Coxeter group $W$ has a unique longest element $w_{0}$.

\begin{definition}
\label{def.coxetersep} Let $\left(  W,S\right)  $ be a Coxeter system with
finite Coxeter group $W$. The \emph{$W$-separable} elements (or, to be
precise, the $\left(  W,S\right)  $-separable elements, since the definition
depends on the whole Coxeter system) are a class of elements of $W$ that is
defined recursively (by induction on $\left\vert S\right\vert $) as follows:
An element $w\in W$ is said to be \emph{$W$-separable} (or \emph{separable
with respect to $W$}) if it

\begin{enumerate}
\item \textbf{either} equals the identity element $1\in W$,

\item \textbf{or} is a $W_{S\setminus\left\{  s\right\}  }$-separable element
of a standard parabolic subgroup $W_{S\setminus\left\{  s\right\}  }$ for some
simple reflection $s\in S$,

\item \textbf{or} has the form $w_{0}u$, where $s \in S$ is some simple
reflection and $u$ is a $W_{S\setminus\left\{  s\right\}  }$-separable element
in the parabolic subgroup $W_{S\setminus\left\{  s\right\}  }$,

\item \textbf{or} has the form $uv$, where $S=I\oplus J$ is an orthogonal
decomposition of $\left(  W,S\right)  $, and where $u\in W_{I}$ and $v\in
W_{J}$ both are separable with respect to their respective parabolic subgroups
(that is, $u$ is $W_{I}$-separable and $v$ is $W_{J}$-separable).
\end{enumerate}

We denote the set of $W$-separable elements by $\operatorname*{Sep}\left(
W,S\right)  $, or by $\operatorname*{Sep}W$ when the simple reflections are
clear from context. We let $\operatorname*{sep}\left(  W,S\right)  =\left\vert
\operatorname*{Sep}\left(  W,S\right)  \right\vert $ denote the number of
$W$-separable elements.
\end{definition}

The terminology \textquotedblleft$W$-separable\textquotedblright\ is slightly
risky, since the notion depends on the whole Coxeter system $\left(
W,S\right)  $ and not just on the abstract group $W$; an example of this
dependence will be given below (Example \ref{exa.D6sep}).

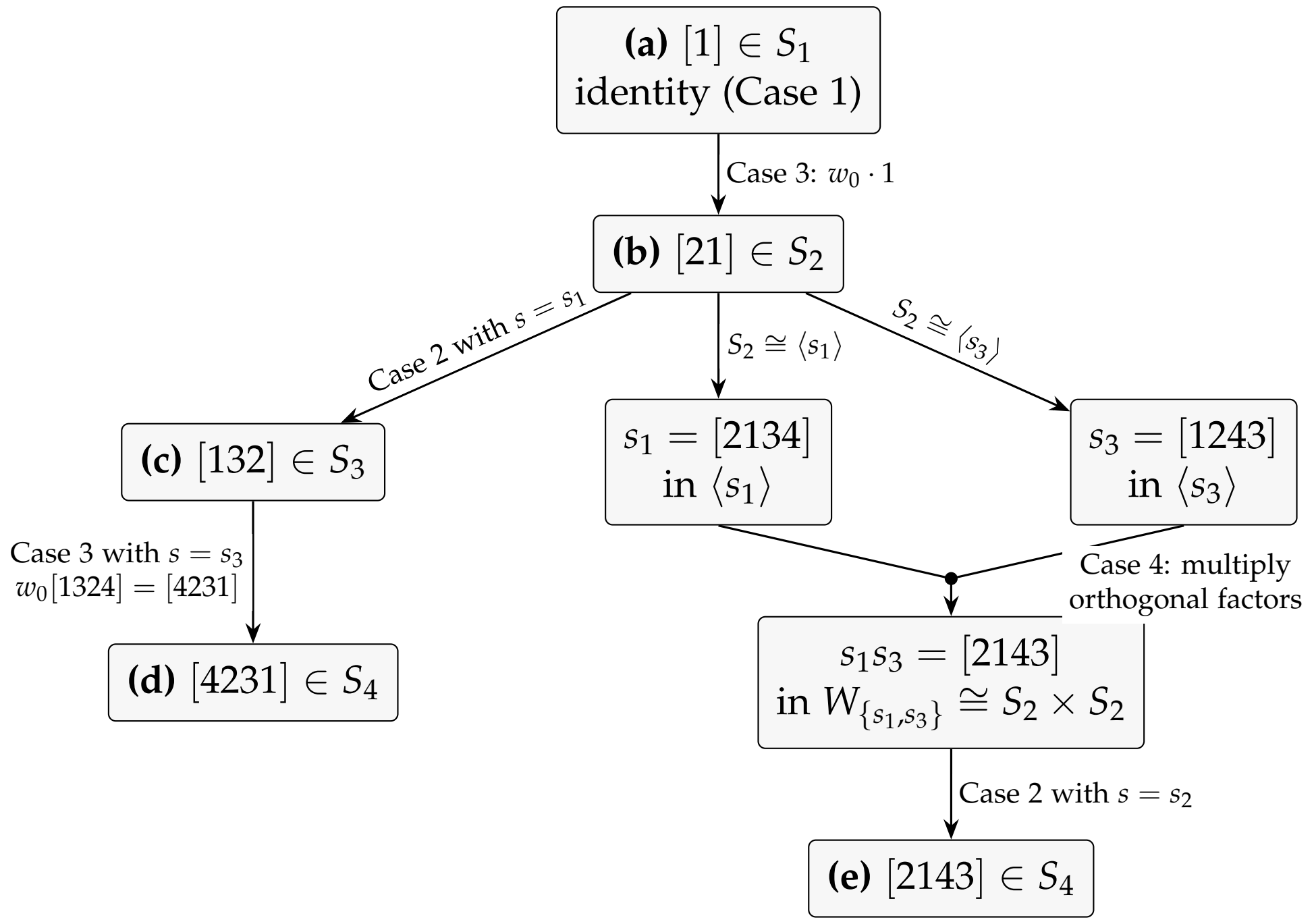
\begin{figure}[h]
\centering
\begin{tikzpicture}[
box/.style={draw,rounded corners=2pt,align=center,inner sep=5pt,
font=\small,fill=black!3},
deduce/.style={-{Stealth[length=2mm]},semithick},
lab/.style={font=\scriptsize,align=center,fill=white,inner sep=2pt},
scale=1.2]
\node[box] (a) at (0,0) {\textbf{(a)} $[1]\in S_1$\\identity (Case 1)};
\node[box] (b) at (0,-1.5) {\textbf{(b)} $[21]\in S_2$};
\node[box] (c) at (-3.8,-3.2) {\textbf{(c)} $[132]\in S_3$};
\node[box] (d) at (-3.8,-5) {\textbf{(d)} $[4231]\in S_4$};
\node[box] (sone) at (0,-3.2) {$s_1=[2134]$\\in $\langle s_1\rangle$};
\node[box] (sthree) at (3.8,-3.2) {$s_3=[1243]$\\in $\langle s_3\rangle$};
\node[box] (product) at (1.9,-5) {$s_1s_3=[2143]$\\in $W_{\{s_1,s_3\}}\cong S_2\times S_2$};
\node[box] (e) at (1.9,-6.6) {\textbf{(e)} $[2143]\in S_4$};
\draw[deduce] (a) -- node[lab,right] {Case 3: $w_0\cdot1$} (b);
\draw[deduce] (b) -- node[lab,sloped,above] {Case 2 with $s=s_1$} (c);
\draw[deduce] (c) -- node[lab,left] {Case 3 with $s=s_3$\\$w_0[1324]=[4231]$} (d);
\draw[deduce] (b) -- node[lab,right] {$S_2\cong\langle s_1\rangle$} (sone);
\draw[deduce] (b) -- node[lab,sloped,above] {$S_2\cong\langle s_3\rangle$} (sthree);
\coordinate (join) at (1.9,-4.15);
\draw[semithick] (sone.south) -- (join) -- (sthree.south);
\fill (join) circle (1.5pt);
\draw[deduce] (join) -- (product.north);
\node[lab,anchor=west] at (2.8,-4.2) {Case 4: multiply\\orthogonal factors};
\draw[deduce] (product) -- node[lab,right] {Case 2 with $s=s_2$} (e);
\end{tikzpicture}
\caption{ The diagram above shows how to build each separable permutation of
Example~\ref{exa.sep.permutations} inductively from the identity permutation.
The case numbers labeling arrows refer to the cases of
Definition~\ref{def.coxetersep}.}%
\label{fig.sep.proof-dag}%
\end{figure}

\begin{example}
\label{exa.sep.permutations}Consider each symmetric group $S_{n}$ as Coxeter
group $A_{n-1}$ with the simple reflections $s_{1},s_{2},\ldots,s_{n-1}$. We
shall write permutations $\sigma\in S_{n}$ in one-line notation, i.e., as
words $\left[  \sigma\left(  1\right)  \ \sigma\left(  2\right)
\ \cdots\ \sigma\left(  n\right)  \right]  $. Then (cf.
Figure~\ref{fig.sep.proof-dag}):

\begin{enumerate}
\item[\textbf{(a)}] The permutation $\left[  1\right]  \in S_{1}$ is $S_{1}%
$-separable, because it equals the identity element $1\in S_{1}$.

\item[\textbf{(b)}] The permutation $\left[  21\right]  \in S_{2}$ is $S_{2}%
$-separable, because it has the form $\left[  21\right]  =w_{0}=w_{0}u$ where
$u=\left[  12\right]  $ is a $\left(  S_{2}\right)  _{S\setminus\left\{
s_{1}\right\}  }$-separable element of the standard parabolic subgroup
$\left(  S_{2}\right)  _{S\setminus\left\{  s_{1}\right\}  }=S_{1}$ (since it
is the identity element of that group).

\item[\textbf{(c)}] The permutation $\left[  132\right]  \in S_{3}$ is $S_{3}%
$-separable, because it is a $\left(  S_{3}\right)  _{S\setminus\left\{
s_{1}\right\}  }$-separable element of the standard parabolic subgroup
$\left(  S_{3}\right)  _{S\setminus\left\{  s_{1}\right\}  }\cong S_{2}$
(indeed, it corresponds to the permutation $\left[  21\right]  \in S_{2}$).

\item[\textbf{(d)}] The permutation $\left[  4231\right]  \in S_{4}$ is
$S_{4}$-separable, because it has the form $w_{0}u$ where $u=\left[
1324\right]  $ is a $\left(  S_{4}\right)  _{S\setminus\left\{  s_{3}\right\}
}$-separable element of the standard parabolic subgroup $\left(  S_{4}\right)
_{S\setminus\left\{  s_{3}\right\}  }\cong S_{3}$ (indeed, it corresponds to
the permutation $\left[  132\right]  \in S_{3}$).

\item[\textbf{(e)}] The permutation $\left[  2143\right]  \in S_{4}$ belongs
to the parabolic subgroup $\left(  S_{4}\right)  _{\left\{  s_{1}\right\}
}\times\left(  S_{4}\right)  _{\left\{  s_{3}\right\}  }\cong S_{2}\times
S_{2}$ of permutations that fix the subsets $\left\{  1,2\right\}  $ and
$\left\{  3,4\right\}  $. As an element of the latter, it is $S_{2}\times
S_{2}$-separable by Case 4 in the definition of $W$-separability, since it
factors as $s_{1}s_{3}$ with $s_{1}\in\left(  S_{4}\right)  _{\left\{
s_{1}\right\}  }$ and $s_{3}\in\left(  S_{4}\right)  _{\left\{  s_{3}\right\}
}$ being separable with respect to their respective Coxeter subgroups. Thus,
it is also $S_{4}$-separable, by Case 2 in the definition of $W$-separability.
\end{enumerate}
\end{example}

\begin{example}
\label{exa.D6sep}Consider the dihedral group with $12$ elements. There are two
ways to realize this group as a Coxeter system. The first is as $G_{2}$, where
$S=\{s_{1},s_{2}\}$ and $m(s_{1},s_{2})=6$. For this Coxeter system, the only
separable elements are $1,s_{1},s_{2},w_{0}s_{2},w_{0}s_{1},w_{0}$. Whereas if
instead we consider the group as $A_{1}\times A_{2}$, that is where
$S=\{s_{1},s_{2},s_{3}\}$ with $m(s_{1},s_{2})=m(s_{1},s_{3})=2$ and
$m(s_{2},s_{3})=3$, then every element of this group is separable.
\end{example}

However, of course, an isomorphism $\left(  W,S\right)  \rightarrow\left(
W^{\prime},S^{\prime}\right)  $ of Coxeter systems does map the $W$-separable
elements to the $W^{\prime}$-separable elements.

The notion of $W$-separable elements is reasonably well-behaved; in
particular, we shall show the following two facts:

\begin{proposition}
\label{prop.w0-inv}Let $\left(  W,S\right)  $ be a Coxeter system with finite
Coxeter group $W$. Let $p\in W$ be $W$-separable. Then, $w_{0}p$ is also $W$-separable.
\end{proposition}

\begin{proposition}
\label{prop.canonical} Let $(W,S)$ be a Coxeter system, where $W$ is finite.
Let $I$ be a subset of $S$. Let $w\in W_{I}$. Then $w$ is $W$-separable if and
only if $w$ is $W_{I}$-separable.
\end{proposition}

The proofs of these (and following) facts will be given in Section
\ref{sec.sep}.

As the name suggests, our notion of $W$-separable elements recovers the
classical separable permutations (see \cite[\textquotedblleft Separable
permutations\textquotedblright]{Ghys16}) in the case when $(W,S)$ is the
type-A Coxeter system $A_{n}$:

\begin{theorem}
\label{thm.separable_equivalence} An element $u$ of $A_{n}$ is $A_{n}%
$-separable as in Definition~\ref{def.coxetersep} if and only if $u$ is
separable as a permutation in $S_{n+1}$.
\end{theorem}

And, generalizing the statement in the introduction, we shall show our main
theorem on separable elements:

\begin{theorem}
\label{thm.sep}Let $\left(  W,S\right)  $ be a Coxeter system with finite
Coxeter group $W$. Let $w\in W$ be $W$-separable. Then, $w$ is $W$-length-symmetric.
\end{theorem}

The reader can find proofs of these results in Section \ref{sec.sep}.
Furthermore, we will prove that our notion of $W$-separability generalizes the
one defined by Gaetz and Gao for finite Weyl groups (Subsection
\ref{subsec.sep.gaga}), and we will count separable elements in type B
(Subsection \ref{subsec.sep.count}).

\section{\label{sec.pfs}Proofs: $W$-length-symmetry}

\subsection{Coefficients}

Before we step to the proofs of the above claims, we establish some basic facts.

We begin with some basic properties of the notation from Definition
\ref{def.coeff}:

\begin{lemma}
\label{lem.coeff}Let $W$ be a group. Let $w,g,h\in W$ and $a\in\mathbf{k}%
\left[  W\right]  $. Then,%
\[
\left[  w\right]  \left(  gah\right)  =\left[  g^{-1}wh^{-1}\right]  \left(
a\right)  .
\]

\end{lemma}

\begin{proof}
Both sides of the claimed equality are $\mathbf{k}$-linear in $a$. Thus, we
WLOG assume that $a$ is an element of the basis $W$ of the $\mathbf{k}$-module
$\mathbf{k}\left[  W\right]  $ (since $a$ is clearly a $\mathbf{k}$-linear
combination of these basis elements). In other words, $a\in W$. Hence, $gah\in
W$ as well, so that
\begin{equation}
\left[  w\right]  \left(  gah\right)  =%
\begin{cases}
1, & \text{if }w=gah;\\
0, & \text{otherwise}%
\end{cases}
\label{pf.lem.coeff.1}%
\end{equation}
(because any $p,q\in W$ satisfy $\left[  p\right]  q=%
\begin{cases}
1, & \text{if }p=q;\\
0, & \text{otherwise}%
\end{cases}
$) and similarly%
\begin{equation}
\left[  g^{-1}wh^{-1}\right]  \left(  a\right)  =%
\begin{cases}
1, & \text{if }g^{-1}wh^{-1}=a;\\
0, & \text{otherwise.}%
\end{cases}
\label{pf.lem.coeff.2}%
\end{equation}
However, $w=gah$ holds if and only if $g^{-1}wh^{-1}=a$. Hence, the right-hand
sides of the equalities (\ref{pf.lem.coeff.1}) and (\ref{pf.lem.coeff.2}) are
equal. Thus, their left-hand sides are equal as well. In other words, $\left[
w\right]  \left(  gah\right)  =\left[  g^{-1}wh^{-1}\right]  \left(  a\right)
$. This proves Lemma \ref{lem.coeff}.
\end{proof}

\subsection{Dihedral groups}

Dihedral groups (Coxeter groups with two generators) are simple but form a
crucial building block of Coxeter group theory. The following lemma is likely
known in some form, but we could not find it in the literature.

\begin{lemma}
\label{lem.dih.sti}Let $\left(  W,S\right)  $ be a Coxeter system, where
$S=\left\{  s,t\right\}  $ is a two-element set. (Thus, $W$ is a dihedral
group.) For each $i\in\mathbb{N}$, we let $\left(  st\right)  _{i}$ denote the
product $ststs\cdots$ with $i$ factors (that is, $\left(  st\right)  ^{i/2}$
if $i$ is even, and $\left(  st\right)  ^{\left(  i-1\right)  /2}s$ if $i$ is
odd), and we let $\left(  ts\right)  _{i}$ denote the product $tstst\cdots$
with $i$ factors. Let $B$ be the $\mathbf{k}$-subalgebra of the group algebra
$\mathbf{k}\left[  W\right]  $ generated by $s+t$. Then, each $i\in\mathbb{N}$
satisfies%
\[
\left(  st\right)  _{i}+\left(  ts\right)  _{i}\in B.
\]

\end{lemma}

\begin{proof}
We prove this by strong induction on $i$. The \textit{base cases} $i=0$ and
$i=1$ are obvious (since $\left(  st\right)  _{0}+\left(  ts\right)
_{0}=1+1=2$ and $\left(  st\right)  _{1}+\left(  ts\right)  _{1}=s+t$). For
the \textit{induction step}, we observe that each $i\geq2$ satisfies%
\begin{align*}
\left(  s+t\right)  \left(  \left(  st\right)  _{i-1}+\left(  ts\right)
_{i-1}\right)   &  =s\underbrace{\left(  st\right)  _{i-1}}_{=s\left(
ts\right)  _{i-2}}+\underbrace{s\left(  ts\right)  _{i-1}}_{=\left(
st\right)  _{i}}+\underbrace{t\left(  st\right)  _{i-1}}_{=\left(  ts\right)
_{i}}+\,t\underbrace{\left(  ts\right)  _{i-1}}_{=t\left(  st\right)  _{i-2}%
}\\
&  =\underbrace{ss}_{=s^{2}=1}\left(  ts\right)  _{i-2}+\left(  st\right)
_{i}+\left(  ts\right)  _{i}+\underbrace{tt}_{=t^{2}=1}\left(  st\right)
_{i-2}\\
&  =\left(  ts\right)  _{i-2}+\left(  st\right)  _{i}+\left(  ts\right)
_{i}+\left(  st\right)  _{i-2}\\
&  =\left(  st\right)  _{i}+\left(  ts\right)  _{i}+\left(  st\right)
_{i-2}+\left(  ts\right)  _{i-2}%
\end{align*}
and thus%
\begin{align*}
\left(  st\right)  _{i}+\left(  ts\right)  _{i}  &  =\underbrace{\left(
s+t\right)  }_{\in B}\underbrace{\left(  \left(  st\right)  _{i-1}+\left(
ts\right)  _{i-1}\right)  }_{\substack{\in B\\\text{(by the induction
hypothesis)}}}-\underbrace{\left(  \left(  st\right)  _{i-2}+\left(
ts\right)  _{i-2}\right)  }_{\substack{\in B\\\text{(by the induction
hypothesis)}}}\\
&  \in B\ \ \ \ \ \ \ \ \ \ \left(  \text{since }B\text{ is a subalgebra}%
\right)  .
\end{align*}
This completes the induction step, and thus the proof of the lemma.
\end{proof}

\begin{proof}
[Proof of Theorem \ref{thm.comm.dihe}.]Let $a,b\in\mathbf{k}$. We must prove
that $\mathbf{L}_{a}\mathbf{L}_{b}=\mathbf{L}_{b}\mathbf{L}_{a}$.

If $\left\vert S\right\vert \leq1$, then the $\mathbf{k}$-algebra
$\mathbf{k}\left[  W\right]  $ is commutative (since it is generated by $S$),
and thus the claim $\mathbf{L}_{a}\mathbf{L}_{b}=\mathbf{L}_{b}\mathbf{L}_{a}$
holds for obvious reasons. Thus, we WLOG assume that $\left\vert S\right\vert
>1$ from now on.

Hence, $\left\vert S\right\vert =2$ (since we assumed that $\left\vert
S\right\vert \leq2$). Write this $2$-element set $S$ as $S=\left\{
s,t\right\}  $.

Use the notations from Lemma \ref{lem.dih.sti}. In particular, define the
$\mathbf{k}$-subalgebra $B$ of $\mathbf{k}\left[  W\right]  $ as in Lemma
\ref{lem.dih.sti}. Note that this algebra $B$ is commutative, since it is
generated by a single element.

However, $W$ is a dihedral group, so it has a well-known description: If we
let $m$ denote the order of $st$ in the group $W$, then the elements of $W$
are precisely the products $\left(  st\right)  _{i}$ and $\left(  ts\right)
_{i}$ for all $i\in\left\{  0,1,\ldots,m\right\}  $, which are furthermore all
distinct except that $\left(  st\right)  _{0}=\left(  ts\right)  _{0}$ and
$\left(  st\right)  _{m}=\left(  ts\right)  _{m}$. Moreover, the products
$\left(  st\right)  _{i}$ and $\left(  ts\right)  _{i}$ have length $i$ for
all $i\in\left\{  0,1,\ldots,m\right\}  $. Altogether, we see that the
elements of $W$ (listed without repetitions) are%
\[
\underbrace{1}_{\text{length }0},\ \underbrace{\left(  st\right)
_{1},\ \left(  ts\right)  _{1}}_{\text{length }1},\ \underbrace{\left(
st\right)  _{2},\ \left(  ts\right)  _{2}}_{\text{length }2},\ \ldots
,\ \underbrace{\left(  st\right)  _{m-1},\ \left(  ts\right)  _{m-1}%
}_{\text{length }m-1},\ \underbrace{\left(  st\right)  _{m}}_{\text{length }%
m}.
\]
Thus, the definition $\mathbf{L}_{a}=\sum_{w\in W}a^{\ell\left(  w\right)  }w$
of $\mathbf{L}_{a}$ can be rewritten as%
\begin{align}
\mathbf{L}_{a}  &  =1+\sum_{i=1}^{m-1}\left(  a^{i}\left(  st\right)
_{i}+a^{i}\left(  ts\right)  _{i}\right)  +a^{m}\left(  st\right)
_{m}\nonumber\\
&  =1+\sum_{i=1}^{m-1}a^{i}\left(  \left(  st\right)  _{i}+\left(  ts\right)
_{i}\right)  +a^{m}\left(  st\right)  _{m}. \label{eq.thm.comm.dihe.1}%
\end{align}
Let us WLOG assume that $\mathbf{k}$ is a polynomial ring over $\mathbb{Q}$
(we can do this, since we are proving a polynomial identity, so we can assume
that $\mathbf{k}$, $a$ and $b$ are generic, as explained in Remark
\ref{rmk.Wls.generic}). Thus, $2$ is invertible in $\mathbf{k}$. Hence, in
(\ref{eq.thm.comm.dihe.1}), we can rewrite $\left(  st\right)  _{m}$ as
$\dfrac{\left(  st\right)  _{m}+\left(  ts\right)  _{m}}{2}$ (since $\left(
st\right)  _{m}=\left(  ts\right)  _{m}$), and thus (\ref{eq.thm.comm.dihe.1})
becomes%
\[
\mathbf{L}_{a}=1+\sum_{i=1}^{m-1}a^{i}\left(  \left(  st\right)  _{i}+\left(
ts\right)  _{i}\right)  +a^{m}\dfrac{\left(  st\right)  _{m}+\left(
ts\right)  _{m}}{2}.
\]
This shows that $\mathbf{L}_{a}$ belongs to $B$ (since Lemma \ref{lem.dih.sti}
shows that all the sums $\left(  st\right)  _{i}+\left(  ts\right)  _{i}$
belong to $B$). Similarly, $\mathbf{L}_{b}$ belongs to $B$ as well. Hence,
$\mathbf{L}_{a}$ and $\mathbf{L}_{b}$ commute (since $B$ is commutative). That
is, $\mathbf{L}_{a}\mathbf{L}_{b}=\mathbf{L}_{b}\mathbf{L}_{a}$. Theorem
\ref{thm.comm.dihe} is thus proved.
\end{proof}

\begin{proof}
[Proof of Corollary \ref{cor.Wls.dih}.]Let $w\in W$. For all $a,b\in
\mathbf{k}$, we have $\mathbf{L}_{a}\mathbf{L}_{b}=\mathbf{L}_{b}%
\mathbf{L}_{a}$ by Theorem \ref{thm.comm.dihe}, and therefore $\left[
w\right]  \left(  \mathbf{L}_{a}\mathbf{L}_{b}\right)  =\left[  w\right]
\left(  \mathbf{L}_{b}\mathbf{L}_{a}\right)  $. In other words, $w$ is
$W$-length-symmetric. This proves Corollary \ref{cor.Wls.dih}.
\end{proof}

\subsection{The antipode, inverses and involutions}

For any group $G$, there is a $\mathbf{k}$-linear map $\mathbf{k}\left[
G\right]  \rightarrow\mathbf{k}\left[  G\right]  $ that sends each group
element $g\in G$ to $g^{-1}$. This map is called the \emph{antipode} of the
group algebra $\mathbf{k}\left[  G\right]  $, and will be denoted by $a\mapsto
a^{\ast}$. It is an involution (i.e., each $a\in\mathbf{k}\left[  G\right]  $
satisfies $a^{\ast\ast}=a$), and is a $\mathbf{k}$-algebra anti-morphism
(i.e., each $a,b\in\mathbf{k}\left[  G\right]  $ satisfy $\left(  ab\right)
^{\ast}=b^{\ast}a^{\ast}$). Thus, it is a $\mathbf{k}$-algebra
anti-isomorphism. Clearly, all $w\in G$ and $a\in\mathbf{k}\left[  G\right]  $
satisfy%
\begin{equation}
\left[  w\right]  \left(  a\right)  =\left[  w^{-1}\right]  \left(  a^{\ast
}\right)  . \label{eq.antipode.coeff}%
\end{equation}

\begin{lemma}
\label{lem.W-S}Let $\left(  W,S\right)  $ be a Coxeter system, where $W$ is
finite. Let $q\in\mathbf{k}$. Then, $\mathbf{L}_{q}^{\ast}=\mathbf{L}_{q}$.
\end{lemma}

\begin{proof}
This is clear, since each $w\in W$ satisfies $\ell\left(  w^{-1}\right)
=\ell\left(  w\right)  $.
\end{proof}

\begin{lemma}
\label{lem.winv}Let $\left(  W,S\right)  $ be a Coxeter system, where $W$ is
finite. Let $a,b\in\mathbf{k}$ and $w\in W$. Then, $\left[  w\right]  \left(
\mathbf{L}_{a}\mathbf{L}_{b}\right)  =\left[  w^{-1}\right]  \left(
\mathbf{L}_{b}\mathbf{L}_{a}\right)  $.
\end{lemma}

\begin{proof}
We have $\left(  \mathbf{L}_{a}\mathbf{L}_{b}\right)  ^{\ast}=\mathbf{L}%
_{b}^{\ast}\mathbf{L}_{a}^{\ast}$ (since the antipode is a $\mathbf{k}%
$-algebra anti-morphism). But Lemma \ref{lem.W-S} yields $\mathbf{L}_{a}%
^{\ast}=\mathbf{L}_{a}$ and $\mathbf{L}_{b}^{\ast}=\mathbf{L}_{b}$. Hence,
$\left(  \mathbf{L}_{a}\mathbf{L}_{b}\right)  ^{\ast}=\underbrace{\mathbf{L}%
_{b}^{\ast}}_{=\mathbf{L}_{b}}\underbrace{\mathbf{L}_{a}^{\ast}}%
_{=\mathbf{L}_{a}}=\mathbf{L}_{b}\mathbf{L}_{a}$. Now,
(\ref{eq.antipode.coeff}) yields%
\[
\left[  w\right]  \left(  \mathbf{L}_{a}\mathbf{L}_{b}\right)  =\left[
w^{-1}\right]  \left(  \left(  \mathbf{L}_{a}\mathbf{L}_{b}\right)  ^{\ast
}\right)  =\left[  w^{-1}\right]  \left(  \mathbf{L}_{b}\mathbf{L}_{a}\right)
\ \ \ \ \ \ \ \ \ \ \left(  \text{since }\left(  \mathbf{L}_{a}\mathbf{L}%
_{b}\right)  ^{\ast}=\mathbf{L}_{b}\mathbf{L}_{a}\right)  ,
\]
and thus Lemma \ref{lem.winv} is proved.
\end{proof}

\begin{proof}
[Proof of Proposition \ref{prop.Wls.inverse}.]Let $a,b\in\mathbf{k}$ be
arbitrary. Then,%
\begin{align*}
\left[  w\right]  \left(  \mathbf{L}_{a}\mathbf{L}_{b}\right)   &  =\left[
w^{-1}\right]  \left(  \mathbf{L}_{b}\mathbf{L}_{a}\right)
\ \ \ \ \ \ \ \ \ \ \left(  \text{by Lemma \ref{lem.winv}}\right)
,\ \ \ \ \ \ \ \ \ \ \text{and}\\
\left[  w\right]  \left(  \mathbf{L}_{b}\mathbf{L}_{a}\right)   &  =\left[
w^{-1}\right]  \left(  \mathbf{L}_{a}\mathbf{L}_{b}\right)
\ \ \ \ \ \ \ \ \ \ \left(  \text{by Lemma \ref{lem.winv}, with the roles of
}a\text{ and }b\text{ swapped}\right)  .
\end{align*}
Hence, the equality $\left[  w\right]  \left(  \mathbf{L}_{a}\mathbf{L}%
_{b}\right)  =\left[  w\right]  \left(  \mathbf{L}_{b}\mathbf{L}_{a}\right)  $
holds if and only if the equality \newline$\left[  w^{-1}\right]  \left(
\mathbf{L}_{b}\mathbf{L}_{a}\right)  =\left[  w^{-1}\right]  \left(
\mathbf{L}_{a}\mathbf{L}_{b}\right)  $ holds. But the former equality (stated
for all $\mathbf{k}$ and $a,b\in\mathbf{k}$) says that $w$ is $W$%
-length-symmetric, whereas the latter equality says that $w^{-1}$ is
$W$-length-symmetric. Hence, $w$ is $W$-length-symmetric if and only if
$w^{-1}$ is $W$-length-symmetric. This proves Proposition
\ref{prop.Wls.inverse}.
\end{proof}

\begin{proof}
[Proof of Proposition \ref{prop.Wls.inv}.]Let $w\in W$ be an involution, and
$a,b\in\mathbf{k}$ be arbitrary. Then, $w$ is an involution, so that
$w^{-1}=w$. However, Lemma \ref{lem.winv} yields%
\[
\left[  w\right]  \left(  \mathbf{L}_{a}\mathbf{L}_{b}\right)  =\left[
w^{-1}\right]  \left(  \mathbf{L}_{b}\mathbf{L}_{a}\right)  =\left[  w\right]
\left(  \mathbf{L}_{b}\mathbf{L}_{a}\right)  \ \ \ \ \ \ \ \ \ \ \left(
\text{since }w^{-1}=w\right)  .
\]
This shows that $w$ is $W$-length-symmetric. This proves Proposition
\ref{prop.Wls.inv}.
\end{proof}

\begin{proof}
[Proof of Proposition \ref{prop.Wls.twinv}.]Forget that we fixed $w$. The
group automorphism $\phi$ of $W$ induces (by linearity) a $\mathbf{k}$-algebra
automorphism $\mathbf{k}\left[  \phi\right]  $ of $\mathbf{k}\left[  W\right]
$, which sends each $w\in W$ to $\phi\left(  w\right)  $.

Moreover, the automorphism $\phi$ of $W$ is a Coxeter group automorphism.
Thus, it preserves the length of any element $w\in W$; that is, we have%
\begin{equation}
\ell\left(  \phi\left(  w\right)  \right)  =\ell\left(  w\right)
\ \ \ \ \ \ \ \ \ \ \text{for each }w\in W. \label{pf.prop.WIs.twinv.pres}%
\end{equation}
Thus, for each $q\in\mathbf{k}$, we have%
\begin{align*}
\left(  \mathbf{k}\left[  \phi\right]  \right)  \left(  \mathbf{L}_{q}\right)
&  =\left(  \mathbf{k}\left[  \phi\right]  \right)  \left(  \sum_{w\in
W}q^{\ell\left(  w\right)  }w\right)  \ \ \ \ \ \ \ \ \ \ \left(  \text{since
}\mathbf{L}_{q}=\sum_{w\in W}q^{\ell\left(  w\right)  }w\right) \\
&  =\sum_{w\in W}\underbrace{q^{\ell\left(  w\right)  }}_{\substack{=q^{\ell
\left(  \phi\left(  w\right)  \right)  }\\\text{(since
(\ref{pf.prop.WIs.twinv.pres}) yields }\ell\left(  w\right)  =\ell\left(
\phi\left(  w\right)  \right)  }}\phi\left(  w\right)  =\sum_{w\in W}%
q^{\ell\left(  \phi\left(  w\right)  \right)  }\phi\left(  w\right) \\
&  =\sum_{w\in W}q^{\ell\left(  w\right)  }w\ \ \ \ \ \ \ \ \ \ \left(
\begin{array}
[c]{c}%
\text{here, we have substituted }w\text{ for }\phi\left(  w\right)  \text{
in}\\
\text{the sum, since }\phi:W\rightarrow W\text{ is a bijection}%
\end{array}
\right) \\
&  =\mathbf{L}_{q}.
\end{align*}
In other words, $\mathbf{k}\left[  \phi\right]  $ fixes $\mathbf{L}_{q}$ for
each $q\in\mathbf{k}$.

Let $a,b\in\mathbf{k}$ be arbitrary. Thus, the $\mathbf{k}$-algebra
automorphism $\mathbf{k}\left[  \phi\right]  $ of $\mathbf{k}\left[  W\right]
$ fixes the elements $\mathbf{L}_{a}$ and $\mathbf{L}_{b}$ (since it fixes
$\mathbf{L}_{q}$ for each $q\in\mathbf{k}$). Being a $\mathbf{k}$-algebra
automorphism, it therefore also fixes the product $\mathbf{L}_{b}%
\mathbf{L}_{a}$. In other words,%
\begin{equation}
\left(  \mathbf{k}\left[  \phi\right]  \right)  \left(  \mathbf{L}%
_{b}\mathbf{L}_{a}\right)  =\mathbf{L}_{b}\mathbf{L}_{a}.
\label{pf.prop.Wls.twinv.0}%
\end{equation}
Moreover, by its very definition, it satisfies%
\begin{equation}
\left[  \phi\left(  w\right)  \right]  \left(  \left(  \mathbf{k}\left[
\phi\right]  \right)  \left(  \alpha\right)  \right)  =\left[  w\right]
\left(  \alpha\right)  \label{pf.prop.Wls.twinv.1}%
\end{equation}
for all $w\in W$ and $\alpha\in\mathbf{k}\left[  W\right]  $.

Now, let $w\in W$ be such that $w^{-1}=\phi\left(  w\right)  $. Now, Lemma
\ref{lem.winv} yields%
\begin{align*}
\left[  w\right]  \left(  \mathbf{L}_{a}\mathbf{L}_{b}\right)   &  =\left[
w^{-1}\right]  \left(  \mathbf{L}_{b}\mathbf{L}_{a}\right) \\
&  =\left[  \phi\left(  w\right)  \right]  \underbrace{\left(  \mathbf{L}%
_{b}\mathbf{L}_{a}\right)  }_{\substack{=\left(  \mathbf{k}\left[
\phi\right]  \right)  \left(  \mathbf{L}_{b}\mathbf{L}_{a}\right)  \\\text{(by
(\ref{pf.prop.Wls.twinv.0}))}}}\ \ \ \ \ \ \ \ \ \ \left(  \text{since }%
w^{-1}=\phi\left(  w\right)  \right) \\
&  =\left[  \phi\left(  w\right)  \right]  \left(  \left(  \mathbf{k}\left[
\phi\right]  \right)  \left(  \mathbf{L}_{b}\mathbf{L}_{a}\right)  \right)
=\left[  w\right]  \left(  \mathbf{L}_{b}\mathbf{L}_{a}\right)
\ \ \ \ \ \ \ \ \ \ \left(  \text{by (\ref{pf.prop.Wls.twinv.1})}\right)  .
\end{align*}
This shows that $w$ is $W$-length-symmetric. This proves Proposition
\ref{prop.Wls.twinv}.
\end{proof}

\subsection{Cosets and parabolic subgroups}

Our next lemma is concerned with cosets and their representatives. Due to the
perennial vexatiousness of their chirality, we recall their definitions:

\begin{definition}
Let $H$ be a subgroup of a group $G$.

\begin{itemize}
\item The \emph{left cosets} of $H$ in $G$ are the subsets of the form $uH$
for $u\in G$. The set of these left cosets is denoted by $G/H$.

A \emph{left transversal} of $H$ in $G$ means a subset $L$ of $G$ that
intersects each of these left cosets in exactly one element (i.e., satisfies
$\left\vert uH\cap L\right\vert =1$ for each $u\in G$). Equivalently, it is a
subset $L$ of $G$ with the property that each $g\in G$ can be written uniquely
as $xh$ with $h\in H$ and $x\in L$.

\item The \emph{right cosets} of $H$ in $G$ are the subsets of the form $Hu$
for $u\in G$. The set of these right cosets is denoted by $H\backslash G$.

A \emph{right transversal} of $H$ in $G$ means a subset $R$ of $G$ that
intersects each of these right cosets in exactly one element (i.e., satisfies
$\left\vert Hu\cap R\right\vert =1$ for each $u\in G$). Equivalently, it is a
subset $R$ of $G$ with the property that each $g\in G$ can be written uniquely
as $hx$ with $h\in H$ and $x\in R$.
\end{itemize}
\end{definition}

The left and right versions here are intimately related: If $L$ is a left
transversal of $H$ in $G$, then $\left\{  x^{-1}\ \mid\ x\in L\right\}  $ is a
right transversal of $H$ in $G$. An analogous statement holds with
\textquotedblleft left\textquotedblright\ and \textquotedblleft
right\textquotedblright\ interchanged.

\begin{lemma}
\label{lem.coset}Let $G$ be a finite group, and let $H$ be a subgroup of $G$.
Let $R$ be a right transversal of $H$ in $G$. For each $x\in R$, let
$p_{x},q_{x}\in\mathbf{k}$ be two scalars. Set%
\[
p=\sum_{x\in R}p_{x}x\ \ \ \ \ \ \ \ \ \ \text{and}\ \ \ \ \ \ \ \ \ \ q=\sum
_{x\in R}q_{x}x^{-1}%
\]
in $\mathbf{k}\left[  G\right]  $. Let $w\in H$ and $\alpha,\beta\in
\mathbf{k}\left[  H\right]  $ (where $\mathbf{k}\left[  H\right]  $ is
embedded in $\mathbf{k}\left[  G\right]  $ in the obvious way). Then,%
\begin{equation}
\left[  w\right]  \left(  \alpha pq\beta\right)  =\left(  \sum_{x\in R}%
p_{x}q_{x}\right)  \left[  w\right]  \left(  \alpha\beta\right)  .
\label{eq.lem.coset.eq}%
\end{equation}

\end{lemma}

\begin{proof}
The claim is $\mathbf{k}$-linear in each of $\alpha$ and $\beta$. Thus, we
WLOG assume that $\alpha,\beta\in H$ (since $\alpha,\beta\in\mathbf{k}\left[
H\right]  $ are clearly $\mathbf{k}$-linear combinations of elements of $H$).
Furthermore, let us set $u:=\alpha^{-1}w\beta^{-1}\in H$ (since $w,\alpha
,\beta\in H$).

Hence, $w,\alpha,\beta\in H\subseteq G$. Thus, Lemma \ref{lem.coeff} (applied
to $\alpha$, $\beta$ and $pq$ instead of $g$, $h$ and $a$) yields%
\begin{equation}
\left[  w\right]  \left(  \alpha pq\beta\right)  =\left[  \alpha^{-1}%
w\beta^{-1}\right]  \left(  pq\right)  =\left[  u\right]  \left(  pq\right)
\label{pf.lem.coset.1}%
\end{equation}
(since $\alpha^{-1}w\beta^{-1}=u$). Moreover, Lemma \ref{lem.coeff} (applied
to $\alpha$, $\beta$ and $1$ instead of $g$, $h$ and $a$) yields%
\begin{equation}
\left[  w\right]  \left(  \alpha\beta\right)  =\left[  \alpha^{-1}w\beta
^{-1}\right]  1=\left[  u\right]  1 \label{pf.lem.coset.2}%
\end{equation}
(since $\alpha^{-1}w\beta^{-1}=u$). Using (\ref{pf.lem.coset.1}) and
(\ref{pf.lem.coset.2}), we can rewrite our claim (\ref{eq.lem.coset.eq}) as%
\begin{equation}
\left[  u\right]  \left(  pq\right)  =\left(  \sum_{x\in R}p_{x}q_{x}\right)
\left[  u\right]  \left(  1\right)  . \label{pf.lem.coset.goalrew}%
\end{equation}
So it remains to prove (\ref{pf.lem.coset.goalrew}).

But the definitions of $p$ and $q$ yield%
\begin{align*}
pq  &  =\left(  \sum_{x\in R}p_{x}x\right)  \left(  \sum_{x\in R}q_{x}%
x^{-1}\right)  =\sum_{\left(  x,y\right)  \in R\times R}p_{x}q_{y}xy^{-1}\\
&  =\sum_{\substack{\left(  x,y\right)  \in R\times R;\\x=y}}p_{x}q_{y}%
xy^{-1}+\sum_{\substack{\left(  x,y\right)  \in R\times R;\\x\neq y}%
}p_{x}q_{y}xy^{-1}%
\end{align*}
(here, we have split the sum into the part with $x=y$ and the part with $x\neq
y$). In view of%
\[
\sum_{\substack{\left(  x,y\right)  \in R\times R;\\x=y}}p_{x}q_{y}%
xy^{-1}=\sum_{x\in R}p_{x}q_{x}\underbrace{xx^{-1}}_{=1}=\sum_{x\in R}%
p_{x}q_{x}1,
\]
we can rewrite this as%
\begin{equation}
pq=\sum_{x\in R}p_{x}q_{x}1+\sum_{\substack{\left(  x,y\right)  \in R\times
R;\\x\neq y}}p_{x}q_{y}xy^{-1}. \label{pf.lem.coset.6}%
\end{equation}
Comparing coefficients of $u$ on both sides of this equality, we obtain%
\begin{equation}
\left[  u\right]  \left(  pq\right)  =\sum_{x\in R}p_{x}q_{x}\left[  u\right]
1+\sum_{\substack{\left(  x,y\right)  \in R\times R;\\x\neq y}}p_{x}%
q_{y}\left[  u\right]  \left(  xy^{-1}\right)  . \label{pf.lem.coset.7}%
\end{equation}

However, we claim that if two elements $x,y\in R$ satisfy $x\neq y$, then
\begin{equation}
\left[  u\right]  \left(  xy^{-1}\right)  =0. \label{pf.lem.coset.8}%
\end{equation}

[\textit{Proof:} Let $x,y\in R$ be two elements satisfying $x\neq y$. Then,
$x$ and $y$ are two distinct elements of $R$. Since $R$ is a right transversal
of $H$, this entails that the right cosets $Hx$ and $Hy$ are distinct. Thus,
the two elements $u$ and $xy^{-1}$ of $G$ are distinct (because if they were
equal, then we would have $xy^{-1}=u\in H$, so that $x\in Hy$ and therefore
$Hx=Hy$, contradicting the distinctness of $Hx$ and $Hy$). Therefore, $\left[
u\right]  \left(  xy^{-1}\right)  =0$ (because $\left[  g\right]  h=0$ for any
two distinct elements $g$ and $h$ of $G$). This proves (\ref{pf.lem.coset.8}).]

Now, (\ref{pf.lem.coset.7}) becomes%
\begin{align*}
\left[  u\right]  \left(  pq\right)   &  =\sum_{x\in R}p_{x}q_{x}\left[
u\right]  1+\sum_{\substack{\left(  x,y\right)  \in R\times R;\\x\neq y}%
}p_{x}q_{y}\underbrace{\left[  u\right]  \left(  xy^{-1}\right)
}_{\substack{=0\\\text{(by (\ref{pf.lem.coset.8}))}}}=\sum_{x\in R}p_{x}%
q_{x}\left[  u\right]  1\\
&  =\left(  \sum_{x\in R}p_{x}q_{x}\right)  \left[  u\right]  \left(
1\right)  .
\end{align*}
This proves (\ref{pf.lem.coset.goalrew}). Thus, Lemma \ref{lem.coset} is proved.
\end{proof}

\begin{proof}
[Proof of Proposition \ref{prop.Wls.WI}.]Given any $q\in\mathbf{k}$, let us
denote the element $\mathbf{L}_{q}\left(  W\right)  \in\mathbf{k}\left[
W\right]  $ by $\mathbf{L}_{q}$, while the analogous element $\mathbf{L}%
_{q}\left(  W_{I}\right)  \in\mathbf{k}\left[  W_{I}\right]  $ shall be
denoted by $\mathbf{L}_{q}^{\prime}$. Thus, $w$ is $W$-length-symmetric if and
only if all $\mathbf{k}$ and $a,b\in\mathbf{k}$ satisfy $\left[  w\right]
\left(  \mathbf{L}_{a}\mathbf{L}_{b}\right)  =\left[  w\right]  \left(
\mathbf{L}_{b}\mathbf{L}_{a}\right)  $, whereas $w$ is $W_{I}$%
-length-symmetric if and only if all $\mathbf{k}$ and $a,b\in\mathbf{k}$
satisfy $\left[  w\right]  \left(  \mathbf{L}_{a}^{\prime}\mathbf{L}%
_{b}^{\prime}\right)  =\left[  w\right]  \left(  \mathbf{L}_{b}^{\prime
}\mathbf{L}_{a}^{\prime}\right)  $. In both cases, we can restrict ourselves
to the \textquotedblleft generic\textquotedblright\ case, i.e., to the case
when $\mathbf{k}$ is the polynomial ring $\mathbb{Z}\left[  a,b\right]  $ and
when $a$ and $b$ are the two indeterminates (as in the proof of Proposition
\ref{prop.Wls.comb}). Let us do so. Thus, $w$ is $W$-length-symmetric if and
only if $\left[  w\right]  \left(  \mathbf{L}_{a}\mathbf{L}_{b}\right)
=\left[  w\right]  \left(  \mathbf{L}_{b}\mathbf{L}_{a}\right)  $, whereas $w$
is $W_{I}$-length-symmetric if and only if $\left[  w\right]  \left(
\mathbf{L}_{a}^{\prime}\mathbf{L}_{b}^{\prime}\right)  =\left[  w\right]
\left(  \mathbf{L}_{b}^{\prime}\mathbf{L}_{a}^{\prime}\right)  $.

Each right coset of $W_{I}$ in $W$ has a unique representative of minimum
length (by \cite[Lemma 9.7 (a)]{Lusztig}). Let $R$ be the set of all these
representatives. Thus, $R$ is a right transversal of $W_{I}$ in $W$. In other
words, each $g\in W$ can be written uniquely as $hx$ with $h\in W_{I}$ and
$x\in R$. That is, the map $W_{I}\times R\rightarrow W,\ \left(  h,x\right)
\mapsto hx$ is a bijection. Moreover, \cite[Lemma 9.7 (b)]{Lusztig} says that
any $x\in R$ and any $h\in W_{I}$ satisfy%
\begin{equation}
\ell\left(  hx\right)  =\ell\left(  h\right)  +\ell\left(  x\right)
\label{eq.prop.Wls.WI.l=l+l}%
\end{equation}
(since $x\in R$ means that $x$ is a representative of minimum length for some
right coset $W_{I}a$ of $W_{I}$ in $W$). Set%
\[
\mathbf{X}_{a}:=\sum_{x\in R}a^{\ell\left(  x\right)  }%
x\ \ \ \ \ \ \ \ \ \ \text{and}\ \ \ \ \ \ \ \ \ \ \mathbf{X}_{b}:=\sum_{x\in
R}b^{\ell\left(  x\right)  }x.
\]
Then, the definition of $\mathbf{L}_{a}=\mathbf{L}_{a}\left(  W\right)  $
yields%
\begin{align}
\mathbf{L}_{a}  &  =\sum_{w\in W}a^{\ell\left(  w\right)  }w=\underbrace{\sum
_{\left(  h,x\right)  \in W_{I}\times R}}_{=\sum_{h\in W_{I}}\ \ \sum_{x\in
R}}\underbrace{a^{\ell\left(  hx\right)  }}_{\substack{=a^{\ell\left(
h\right)  +\ell\left(  x\right)  }\\\text{(by (\ref{eq.prop.Wls.WI.l=l+l}))}%
}}hx\nonumber\\
&  \ \ \ \ \ \ \ \ \ \ \ \ \ \ \ \ \ \ \ \ \left(
\begin{array}
[c]{c}%
\text{here, we have substituted }hx\text{ for }w\text{ in the sum,}\\
\text{since the map }W_{I}\times R\rightarrow W,\ \left(  h,x\right)  \mapsto
hx\text{ is a bijection}%
\end{array}
\right) \nonumber\\
&  =\sum_{h\in W_{I}}\ \ \sum_{x\in R}a^{\ell\left(  h\right)  +\ell\left(
x\right)  }hx=\underbrace{\left(  \sum_{h\in W_{I}}a^{\ell\left(  h\right)
}h\right)  }_{\substack{=\mathbf{L}_{a}^{\prime}\\\text{(by the definition of
}\mathbf{L}_{a}^{\prime}=\mathbf{L}_{a}\left(  W_{I}\right)  \text{)}%
}}\ \ \underbrace{\left(  \sum_{x\in R}a^{\ell\left(  x\right)  }x\right)
}_{\substack{=\mathbf{X}_{a}\\\text{(by the definition of }\mathbf{X}%
_{a}\text{)}}}\nonumber\\
&  =\mathbf{L}_{a}^{\prime}\mathbf{X}_{a}. \label{pf.prop.Wls.WI.3}%
\end{align}
Similarly, $\mathbf{L}_{b}=\mathbf{L}_{b}^{\prime}\mathbf{X}_{b}$. Applying
the antipode to this latter equality, we obtain%
\begin{equation}
\mathbf{L}_{b}^{\ast}=\left(  \mathbf{L}_{b}^{\prime}\mathbf{X}_{b}\right)
^{\ast}=\mathbf{X}_{b}^{\ast}\mathbf{L}_{b}^{\prime\ast}
\label{pf.prop.Wls.WI.4}%
\end{equation}
(since the antipode is a $\mathbf{k}$-algebra anti-morphism). But Lemma
\ref{lem.W-S} yields $\mathbf{L}_{b}^{\ast}=\mathbf{L}_{b}$ and $\mathbf{L}%
_{b}^{\prime\ast}=\mathbf{L}_{b}^{\prime}$. Hence, we can rewrite
(\ref{pf.prop.Wls.WI.4}) as
\begin{equation}
\mathbf{L}_{b}=\mathbf{X}_{b}^{\ast}\mathbf{L}_{b}^{\prime}.
\label{pf.prop.Wls.WI.4b}%
\end{equation}

However, $w\in W_{I}$ and $\mathbf{L}_{a}^{\prime},\mathbf{L}_{b}^{\prime}%
\in\mathbf{k}\left[  W_{I}\right]  $ and
\[
\mathbf{X}_{a}=\sum_{x\in R}a^{\ell\left(  x\right)  }%
x\ \ \ \ \ \ \ \ \ \ \text{and}\ \ \ \ \ \ \ \ \ \ \mathbf{X}_{b}^{\ast}%
=\sum_{x\in R}b^{\ell\left(  x\right)  }x^{-1}%
\]
(the latter equality follows from $\mathbf{X}_{b}=\sum_{x\in R}b^{\ell\left(
x\right)  }x$ by the definition of the antipode). Hence, Lemma \ref{lem.coset}
(applied to $G=W$ and $H=W_{I}$ and $p_{x}=a^{\ell\left(  x\right)  }$ and
$q_{x}=b^{\ell\left(  x\right)  }$ and $p=\mathbf{X}_{a}$ and $q=\mathbf{X}%
_{b}^{\ast}$ and $\alpha=\mathbf{L}_{a}^{\prime}$ and $\beta=\mathbf{L}%
_{b}^{\prime}$) yields%
\[
\left[  w\right]  \left(  \mathbf{L}_{a}^{\prime}\mathbf{X}_{a}\mathbf{X}%
_{b}^{\ast}\mathbf{L}_{b}^{\prime}\right)  =\left(  \sum_{x\in R}%
a^{\ell\left(  x\right)  }b^{\ell\left(  x\right)  }\right)  \left[  w\right]
\left(  \mathbf{L}_{a}^{\prime}\mathbf{L}_{b}^{\prime}\right)  .
\]
Using (\ref{pf.prop.Wls.WI.3}) and (\ref{pf.prop.Wls.WI.4b}), we can rewrite
this as%
\begin{equation}
\left[  w\right]  \left(  \mathbf{L}_{a}\mathbf{L}_{b}\right)  =\left(
\sum_{x\in R}a^{\ell\left(  x\right)  }b^{\ell\left(  x\right)  }\right)
\left[  w\right]  \left(  \mathbf{L}_{a}^{\prime}\mathbf{L}_{b}^{\prime
}\right)  . \label{pf.prop.Wls.WI.6a}%
\end{equation}
The same argument (with the roles of $a$ and $b$ interchanged) shows that%
\begin{equation}
\left[  w\right]  \left(  \mathbf{L}_{b}\mathbf{L}_{a}\right)  =\left(
\sum_{x\in R}b^{\ell\left(  x\right)  }a^{\ell\left(  x\right)  }\right)
\left[  w\right]  \left(  \mathbf{L}_{b}^{\prime}\mathbf{L}_{a}^{\prime
}\right)  . \label{pf.prop.Wls.WI.6b}%
\end{equation}

Clearly, the factors $\sum_{x\in R}a^{\ell\left(  x\right)  }b^{\ell\left(
x\right)  }$ and $\sum_{x\in R}b^{\ell\left(  x\right)  }a^{\ell\left(
x\right)  }$ on the right-hand sides of the two equalities
(\ref{pf.prop.Wls.WI.6a}) and (\ref{pf.prop.Wls.WI.6b}) are equal, and are not
zero divisors (since the polynomial ring $\mathbb{Z}\left[  a,b\right]  $ is
an integral domain, and since $R$ is nonempty, so that $\sum_{x\in R}%
a^{\ell\left(  x\right)  }b^{\ell\left(  x\right)  }\neq0$ in $\mathbb{Z}%
\left[  a,b\right]  $). Hence, these two equalities show that $\left[
w\right]  \left(  \mathbf{L}_{a}\mathbf{L}_{b}\right)  =\left[  w\right]
\left(  \mathbf{L}_{b}\mathbf{L}_{a}\right)  $ if and only if $\left[
w\right]  \left(  \mathbf{L}_{a}^{\prime}\mathbf{L}_{b}^{\prime}\right)
=\left[  w\right]  \left(  \mathbf{L}_{b}^{\prime}\mathbf{L}_{a}^{\prime
}\right)  $.

In other words, $w$ is $W$-length-symmetric if and only if $w$ is $W_{I}%
$-length-symmetric (since we know that $w$ is $W$-length-symmetric if and only
if $\left[  w\right]  \left(  \mathbf{L}_{a}\mathbf{L}_{b}\right)  =\left[
w\right]  \left(  \mathbf{L}_{b}\mathbf{L}_{a}\right)  $, whereas $w$ is
$W_{I}$-length-symmetric if and only if $\left[  w\right]  \left(
\mathbf{L}_{a}^{\prime}\mathbf{L}_{b}^{\prime}\right)  =\left[  w\right]
\left(  \mathbf{L}_{b}^{\prime}\mathbf{L}_{a}^{\prime}\right)  $). This proves
Proposition \ref{prop.Wls.WI}.
\end{proof}

\begin{proof}
[Proof of Proposition \ref{prop.Wls.len2}.]From $\ell\left(  w\right)  \leq2$,
we know that we can write $w$ as a product of at most two simple reflections.
Let $I$ be the set of all simple reflections appearing in this product. Then,
$w\in W_{I}$ and $\left\vert I\right\vert \leq2$. Hence, Corollary
\ref{cor.Wls.dih} (applied to the Coxeter system $\left(  W_{I},I\right)  $
instead of $\left(  W,S\right)  $) yields that $w$ is $W_{I}$%
-length-symmetric. By Proposition \ref{prop.Wls.WI}, this entails that $w$ is
$W$-length-symmetric. Thus, Proposition \ref{prop.Wls.len2} is proved.
\end{proof}

\subsection{Orthogonal decompositions}

\begin{proof}
[Proof of Proposition \ref{prop.Wls.I+J}.]Let $\mathbf{k}$ be any commutative
ring, and $a,b\in\mathbf{k}$.

Since $S=I\oplus J$, we have $W\cong W_{I}\times W_{J}$; more specifically,
the map%
\begin{align}
W_{I}\times W_{J}  &  \rightarrow W,\nonumber\\
\left(  x,y\right)   &  \mapsto xy \label{pf.prop.Wls.I+J.bijgp}%
\end{align}
is a group isomorphism. This isomorphism entails that the $\mathbf{k}$-linear
map%
\begin{align}
\Phi:\mathbf{k}\left[  W_{I}\right]  \otimes\mathbf{k}\left[  W_{J}\right]
&  \rightarrow\mathbf{k}\left[  W\right]  ,\nonumber\\
x\otimes y  &  \mapsto xy \label{pf.prop.Wls.I+J.bijal}%
\end{align}
is a $\mathbf{k}$-algebra isomorphism (because there is a $\mathbf{k}$-algebra
isomorphism%
\begin{align*}
\mathbf{k}\left[  W_{I}\times W_{J}\right]   &  \rightarrow\mathbf{k}\left[
W_{I}\right]  \otimes\mathbf{k}\left[  W_{J}\right]  ,\\
\left(  x,y\right)   &  \mapsto x\otimes y,
\end{align*}
and if we identify $\mathbf{k}\left[  W_{I}\right]  \otimes\mathbf{k}\left[
W_{J}\right]  $ with $\mathbf{k}\left[  W_{I}\times W_{J}\right]  $ along this
isomorphism, then (\ref{pf.prop.Wls.I+J.bijal}) is the linearization of the
group isomorphism (\ref{pf.prop.Wls.I+J.bijgp})). The construction of $\Phi$
shows that any $\alpha\in\mathbf{k}\left[  W_{I}\right]  $ and $\beta
\in\mathbf{k}\left[  W_{J}\right]  $ satisfy%
\begin{equation}
\left[  uv\right]  \left(  \Phi\left(  \alpha\otimes\beta\right)  \right)
=\left[  u\right]  \left(  \alpha\right)  \cdot\left[  v\right]  \left(
\beta\right)  . \label{pf.prop.Wls.I+J.uvPhi}%
\end{equation}
(Indeed, this equality is $\mathbf{k}$-linear in $\alpha$ and $\beta$, and
thus needs only be checked for all $\alpha\in W_{I}$ and $\beta\in W_{J}$; but
then it boils down to \textquotedblleft$\alpha\beta=uv$ if and only if
$\alpha=u$ and $\beta=v$\textquotedblright.)

Moreover, any element of $W_{I}$ commutes with any element of $W_{J}$; thus,
for any $x\in W_{I}$ and $y\in W_{J}$, we have%
\begin{equation}
\ell\left(  xy\right)  =\ell\left(  x\right)  +\ell\left(  y\right)
\label{pf.prop.Wls.I+J.l=l+l}%
\end{equation}
(since any reduced expression for $xy$ can be transformed into a reduced
expression for $x$ times a reduced expression for $y$ by commuting its
$I$-factors to the left of all its $J$-factors\footnote{In more detail: Let
$x\in W_{I}$ and $y\in W_{J}$. Let $s_{1}s_{2}\cdots s_{p}$ be any reduced
expression for $xy$. Then, $p=\ell\left(  xy\right)  $. Moreover, each $s_{i}$
belongs to $S=I\sqcup J$, hence lies in either $I$ or $J$. Since each $i\in I$
commutes with each $j\in J$ (because $S=I\oplus J$ is an orthogonal
decomposition), we can freely move all the $s_{i}$ that belong to $I$ in the
reduced expression $s_{1}s_{2}\cdots s_{p}$ past the $s_{i}$ that belong to
$J$. In particular, we can move all the former $s_{i}$ to the left end of the
reduced expression, while letting the latter $s_{i}$ accumulate on the right
end. As a consequence, our reduced expression assumes the form $t_{1}%
t_{2}\cdots t_{p}$, where the first $k$ factors $t_{1},t_{2},\ldots,t_{k}$
(for some $k\in\left\{  0,1,\ldots,p\right\}  $) belong to $I$ while the
remaining $p-k$ factors $t_{k+1},t_{k+2},\ldots,t_{p}$ belong to $J$. Thus, of
course, $t_{1}t_{2}\cdots t_{k}\in W_{I}$ and $t_{k+1}t_{k+2}\cdots t_{p}\in
W_{J}$. Hence,%
\[
xy=s_{1}s_{2}\cdots s_{p}=t_{1}t_{2}\cdots t_{p}=\underbrace{\left(
t_{1}t_{2}\cdots t_{k}\right)  }_{\in W_{I}}\underbrace{\left(  t_{k+1}%
t_{k+2}\cdots t_{p}\right)  }_{\in W_{J}}.
\]
In other words, the map (\ref{pf.prop.Wls.I+J.bijgp}) sends the pairs $\left(
x,y\right)  \in W_{I}\times W_{J}$ and $\left(  t_{1}t_{2}\cdots
t_{k},\ t_{k+1}t_{k+2}\cdots t_{p}\right)  \in W_{I}\times W_{J}$ to the same
image. Since this map is injective (being a group isomorphism), it thus
follows that these pairs are equal. That is, $x=t_{1}t_{2}\cdots t_{k}$ and
$y=t_{k+1}t_{k+2}\cdots t_{p}$. Hence, $\ell\left(  x\right)  \leq k$ and
$\ell\left(  y\right)  \leq p-k$. Summing these two inequalities, we obtain
$\ell\left(  x\right)  +\ell\left(  y\right)  \leq k+\left(  p-k\right)
=p=\ell\left(  xy\right)  $. Combining this with $\ell\left(  xy\right)
\leq\ell\left(  x\right)  +\ell\left(  y\right)  $ (which is a general
property of Coxeter lengths in any Coxeter group), we obtain $\ell\left(
xy\right)  =\ell\left(  x\right)  +\ell\left(  y\right)  $. This proves
(\ref{pf.prop.Wls.I+J.l=l+l}).}). But Definition \ref{def.Lq} yields%
\[
\mathbf{L}_{a}\left(  W_{I}\right)  =\sum_{x\in W_{I}}a^{\ell\left(  x\right)
}x\ \ \ \ \ \ \ \ \ \ \text{and}\ \ \ \ \ \ \ \ \ \ \mathbf{L}_{a}\left(
W_{J}\right)  =\sum_{y\in W_{J}}a^{\ell\left(  y\right)  }y,
\]
so that%
\begin{align*}
\mathbf{L}_{a}\left(  W_{I}\right)  \otimes\mathbf{L}_{a}\left(  W_{J}\right)
&  =\left(  \sum_{x\in W_{I}}a^{\ell\left(  x\right)  }x\right)
\otimes\left(  \sum_{y\in W_{J}}a^{\ell\left(  y\right)  }y\right) \\
&  =\sum_{\left(  x,y\right)  \in W_{I}\times W_{J}}\ \ \underbrace{a^{\ell
\left(  x\right)  }a^{\ell\left(  y\right)  }}_{\substack{=a^{\ell\left(
x\right)  +\ell\left(  y\right)  }=a^{\ell\left(  xy\right)  }\\\text{(by
(\ref{pf.prop.Wls.I+J.l=l+l}))}}}x\otimes y=\sum_{\left(  x,y\right)  \in
W_{I}\times W_{J}}a^{\ell\left(  xy\right)  }x\otimes y.
\end{align*}
Applying the linear map $\Phi$ to this equality, we obtain%
\begin{align*}
\Phi\left(  \mathbf{L}_{a}\left(  W_{I}\right)  \otimes\mathbf{L}_{a}\left(
W_{J}\right)  \right)   &  =\sum_{\left(  x,y\right)  \in W_{I}\times W_{J}%
}a^{\ell\left(  xy\right)  }\underbrace{\Phi\left(  x\otimes y\right)  }%
_{=xy}=\sum_{\left(  x,y\right)  \in W_{I}\times W_{J}}a^{\ell\left(
xy\right)  }xy\\
&  =\sum_{w\in W}a^{\ell\left(  w\right)  }w\ \ \ \ \ \ \ \ \ \ \left(
\begin{array}
[c]{c}%
\text{here, we have substituted }w\text{ for }xy\text{,}\\
\text{since the map (\ref{pf.prop.Wls.I+J.bijgp}) is a bijection}%
\end{array}
\right) \\
&  =\mathbf{L}_{a}\left(  W\right)  .
\end{align*}
Similarly, $\Phi\left(  \mathbf{L}_{b}\left(  W_{I}\right)  \otimes
\mathbf{L}_{b}\left(  W_{J}\right)  \right)  =\mathbf{L}_{b}\left(  W\right)
$. Since $\Phi$ is an algebra morphism, we now have%
\begin{align*}
&  \Phi\left(  \left(  \mathbf{L}_{a}\left(  W_{I}\right)  \otimes
\mathbf{L}_{a}\left(  W_{J}\right)  \right)  \left(  \mathbf{L}_{b}\left(
W_{I}\right)  \otimes\mathbf{L}_{b}\left(  W_{J}\right)  \right)  \right) \\
&  =\underbrace{\Phi\left(  \mathbf{L}_{a}\left(  W_{I}\right)  \otimes
\mathbf{L}_{a}\left(  W_{J}\right)  \right)  }_{=\mathbf{L}_{a}\left(
W\right)  }\cdot\underbrace{\Phi\left(  \mathbf{L}_{b}\left(  W_{I}\right)
\otimes\mathbf{L}_{b}\left(  W_{J}\right)  \right)  }_{=\mathbf{L}_{b}\left(
W\right)  }\\
&  =\mathbf{L}_{a}\left(  W\right)  \cdot\mathbf{L}_{b}\left(  W\right)  ,
\end{align*}
so that%
\begin{align*}
&  \mathbf{L}_{a}\left(  W\right)  \cdot\mathbf{L}_{b}\left(  W\right) \\
&  =\Phi\left(  \underbrace{\left(  \mathbf{L}_{a}\left(  W_{I}\right)
\otimes\mathbf{L}_{a}\left(  W_{J}\right)  \right)  \left(  \mathbf{L}%
_{b}\left(  W_{I}\right)  \otimes\mathbf{L}_{b}\left(  W_{J}\right)  \right)
}_{=\mathbf{L}_{a}\left(  W_{I}\right)  \cdot\mathbf{L}_{b}\left(
W_{I}\right)  \otimes\mathbf{L}_{a}\left(  W_{J}\right)  \cdot\mathbf{L}%
_{b}\left(  W_{J}\right)  }\right) \\
&  =\Phi\left(  \mathbf{L}_{a}\left(  W_{I}\right)  \cdot\mathbf{L}_{b}\left(
W_{I}\right)  \otimes\mathbf{L}_{a}\left(  W_{J}\right)  \cdot\mathbf{L}%
_{b}\left(  W_{J}\right)  \right)  .
\end{align*}
Thus,%
\begin{align}
&  \left[  uv\right]  \left(  \mathbf{L}_{a}\left(  W\right)  \cdot
\mathbf{L}_{b}\left(  W\right)  \right) \nonumber\\
&  =\left[  uv\right]  \left(  \Phi\left(  \mathbf{L}_{a}\left(  W_{I}\right)
\cdot\mathbf{L}_{b}\left(  W_{I}\right)  \otimes\mathbf{L}_{a}\left(
W_{J}\right)  \cdot\mathbf{L}_{b}\left(  W_{J}\right)  \right)  \right)
\nonumber\\
&  =\left[  u\right]  \left(  \mathbf{L}_{a}\left(  W_{I}\right)
\cdot\mathbf{L}_{b}\left(  W_{I}\right)  \right)  \cdot\left[  v\right]
\left(  \mathbf{L}_{a}\left(  W_{J}\right)  \cdot\mathbf{L}_{b}\left(
W_{J}\right)  \right)  \label{pf.prop.Wls.I+J.6a}%
\end{align}
(by (\ref{pf.prop.Wls.I+J.uvPhi}), applied to $\alpha=\mathbf{L}_{a}\left(
W_{I}\right)  \cdot\mathbf{L}_{b}\left(  W_{I}\right)  $ and $\beta
=\mathbf{L}_{a}\left(  W_{J}\right)  \cdot\mathbf{L}_{b}\left(  W_{J}\right)
$). The same argument (with the roles of $a$ and $b$ interchanged) shows that%
\begin{align}
&  \left[  uv\right]  \left(  \mathbf{L}_{b}\left(  W\right)  \cdot
\mathbf{L}_{a}\left(  W\right)  \right) \nonumber\\
&  =\left[  u\right]  \left(  \mathbf{L}_{b}\left(  W_{I}\right)
\cdot\mathbf{L}_{a}\left(  W_{I}\right)  \right)  \cdot\left[  v\right]
\left(  \mathbf{L}_{b}\left(  W_{J}\right)  \cdot\mathbf{L}_{a}\left(
W_{J}\right)  \right)  . \label{pf.prop.Wls.I+J.6b}%
\end{align}

However, $\left[  u\right]  \left(  \mathbf{L}_{a}\left(  W_{I}\right)
\cdot\mathbf{L}_{b}\left(  W_{I}\right)  \right)  =\left[  u\right]  \left(
\mathbf{L}_{b}\left(  W_{I}\right)  \cdot\mathbf{L}_{a}\left(  W_{I}\right)
\right)  $ (since $u$ is $W_{I}$-length-symmetric) and $\left[  v\right]
\left(  \mathbf{L}_{a}\left(  W_{J}\right)  \cdot\mathbf{L}_{b}\left(
W_{J}\right)  \right)  =\left[  v\right]  \left(  \mathbf{L}_{b}\left(
W_{J}\right)  \cdot\mathbf{L}_{a}\left(  W_{J}\right)  \right)  $ (since $v$
is $W_{J}$-length-symmetric). Thus, the right-hand sides of
(\ref{pf.prop.Wls.I+J.6a}) and (\ref{pf.prop.Wls.I+J.6b}) are equal (factor by
factor). Consequently, so are the left-hand sides. In other words,
\[
\left[  uv\right]  \left(  \mathbf{L}_{a}\left(  W\right)  \cdot\mathbf{L}%
_{b}\left(  W\right)  \right)  =\left[  uv\right]  \left(  \mathbf{L}%
_{b}\left(  W\right)  \cdot\mathbf{L}_{a}\left(  W\right)  \right)  .
\]
This shows that $uv$ is $W$-length-symmetric. This proves Proposition
\ref{prop.Wls.I+J}.
\end{proof}

\subsection{The longest word}

\begin{lemma}
\label{lem.w0Lq}Let $\left(  W,S\right)  $ be a Coxeter system, where $W$ is
finite. Let $q\in\mathbf{k}$ be invertible. Then,%
\[
w_{0}\mathbf{L}_{q}=\mathbf{L}_{q}w_{0}=q^{\ell\left(  w_{0}\right)
}\mathbf{L}_{1/q}.
\]

\end{lemma}

\begin{proof}
As we know, $w_{0}$ is an involution; thus, $w_{0}w_{0}=1$ and $w_{0}%
^{-1}=w_{0}$. Moreover,%
\begin{equation}
\ell\left(  w_{0}w\right)  =\ell\left(  w_{0}\right)  -\ell\left(  w\right)
\label{pf.lem.w0Lq.1}%
\end{equation}
for any $w\in W$ (see, e.g., \cite[\S 11.2 (b)]{Lusztig} for a proof). But the
definition of $\mathbf{L}_{1/q}$ yields
\begin{equation}
\mathbf{L}_{1/q}=\sum_{w\in W}\left(  1/q\right)  ^{\ell\left(  w\right)  }w.
\label{pf.lem.w0Lq.2}%
\end{equation}
Now, the definition of $\mathbf{L}_{q}$ yields%
\begin{align*}
w_{0}\mathbf{L}_{q}  &  =w_{0}\sum_{w\in W}q^{\ell\left(  w\right)  }%
w=\sum_{w\in W}q^{\ell\left(  w\right)  }w_{0}w\\
&  =\sum_{w\in W}\underbrace{q^{\ell\left(  w_{0}w\right)  }}%
_{\substack{=q^{\ell\left(  w_{0}\right)  -\ell\left(  w\right)  }\\\text{(by
(\ref{pf.lem.w0Lq.1}))}}}\underbrace{w_{0}w_{0}}_{=1}%
w\ \ \ \ \ \ \ \ \ \ \left(
\begin{array}
[c]{c}%
\text{here, we have substituted }w_{0}w\\
\text{for }w\text{ in the sum}%
\end{array}
\right) \\
&  =\sum_{w\in W}\underbrace{q^{\ell\left(  w_{0}\right)  -\ell\left(
w\right)  }}_{\substack{=q^{\ell\left(  w_{0}\right)  }/q^{\ell\left(
w\right)  }\\=q^{\ell\left(  w_{0}\right)  }\left(  1/q\right)  ^{\ell\left(
w\right)  }}}w=q^{\ell\left(  w_{0}\right)  }\underbrace{\sum_{w\in W}\left(
1/q\right)  ^{\ell\left(  w\right)  }w}_{\substack{=\mathbf{L}_{1/q}%
\\\text{(by (\ref{pf.lem.w0Lq.2}))}}}=q^{\ell\left(  w_{0}\right)  }%
\mathbf{L}_{1/q}.
\end{align*}
Applying the antipode to both sides of this equality, we obtain%
\[
\left(  w_{0}\mathbf{L}_{q}\right)  ^{\ast}=q^{\ell\left(  w_{0}\right)
}\mathbf{L}_{1/q}^{\ast}.
\]
Since $\left(  w_{0}\mathbf{L}_{q}\right)  ^{\ast}=\mathbf{L}_{q}^{\ast}%
w_{0}^{\ast}$ (because the antipode is an algebra anti-morphism), we can
rewrite this as%
\[
\mathbf{L}_{q}^{\ast}w_{0}^{\ast}=q^{\ell\left(  w_{0}\right)  }%
\mathbf{L}_{1/q}^{\ast}.
\]
Using Lemma \ref{lem.W-S} and the obvious fact that $w_{0}^{\ast}=w_{0}%
^{-1}=w_{0}$, we can rewrite this as%
\[
\mathbf{L}_{q}w_{0}=q^{\ell\left(  w_{0}\right)  }\mathbf{L}_{1/q}.
\]
Combining this with $w_{0}\mathbf{L}_{q}=q^{\ell\left(  w_{0}\right)
}\mathbf{L}_{1/q}$, we obtain the whole claim of Lemma \ref{lem.w0Lq}.
\end{proof}

\begin{proof}
[Proof of Proposition \ref{prop.Wls.ref}.]Recall that $w_{0}$ is an
involution, so that $w_{0}=w_{0}^{-1}$ and $w_{0}w_{0}=1$. Hence,
$w_{0}\left(  w_{0}w\right)  =\underbrace{w_{0}w_{0}}_{=1}w=w$. Thus, the
relation between $w$ and $w_{0}w$ is mutual. Hence, we only need to prove the
\textquotedblleft only if\textquotedblright\ direction of Proposition
\ref{prop.Wls.ref} (as the \textquotedblleft if\textquotedblright\ direction
will then follow by applying the \textquotedblleft only if\textquotedblright%
\ direction to $w_{0}w$ instead of $w$).

So let us assume that $w$ is $W$-length-symmetric. We must show that $w_{0}w$
is $W$-length-symmetric. Set $u:=w_{0}w$. Thus, we must show that $u$ is $W$-length-symmetric.

Fix a commutative ring $\mathbf{k}$ and $a,b\in\mathbf{k}$. We must show that
$\left[  u\right]  \left(  \mathbf{L}_{a}\mathbf{L}_{b}\right)  =\left[
u\right]  \left(  \mathbf{L}_{b}\mathbf{L}_{a}\right)  $. We WLOG assume that
$\mathbf{k}$ is a ring of Laurent polynomials and $a$ and $b$ are two
indeterminates (since the claim we are proving is a polynomial identity in $a$
and $b$). Thus, $a$ and $b$ are invertible in $\mathbf{k}$. Since $w$ is
$W$-length-symmetric, we have $\left[  w\right]  \left(  \mathbf{L}%
_{a}\mathbf{L}_{b}\right)  =\left[  w\right]  \left(  \mathbf{L}_{b}%
\mathbf{L}_{a}\right)  $ and likewise
\[
\left[  w\right]  \left(  \mathbf{L}_{a}\mathbf{L}_{1/b}\right)  =\left[
w\right]  \left(  \mathbf{L}_{1/b}\mathbf{L}_{a}\right)  .
\]

From Lemma \ref{lem.w0Lq}, we have%
\[
w_{0}\mathbf{L}_{a}=\mathbf{L}_{a}w_{0}=a^{\ell\left(  w_{0}\right)
}\mathbf{L}_{1/a}%
\]
and%
\[
w_{0}\mathbf{L}_{b}=\mathbf{L}_{b}w_{0}=b^{\ell\left(  w_{0}\right)
}\mathbf{L}_{1/b}.
\]

From $u=w_{0}w=w_{0}^{-1}w$ (since $w_{0}=w_{0}^{-1}$), we obtain%
\[
\left[  u\right]  \left(  \mathbf{L}_{a}\mathbf{L}_{b}\right)  =\left[
w_{0}^{-1}w\right]  \left(  \mathbf{L}_{a}\mathbf{L}_{b}\right)  =\left[
w\right]  \left(  w_{0}\mathbf{L}_{a}\mathbf{L}_{b}\right)
\]
(since Lemma \ref{lem.coeff} (applied to $w_{0}$, $1$ and $\mathbf{L}%
_{a}\mathbf{L}_{b}$ instead of $g$, $h$ and $a$) yields $\left[  w\right]
\left(  w_{0}\mathbf{L}_{a}\mathbf{L}_{b}\right)  =\left[  w_{0}^{-1}w\right]
\left(  \mathbf{L}_{a}\mathbf{L}_{b}\right)  $). Therefore,%
\[
\left[  u\right]  \left(  \mathbf{L}_{a}\mathbf{L}_{b}\right)  =\left[
w\right]  \left(  \underbrace{w_{0}\mathbf{L}_{a}}_{=\mathbf{L}_{a}w_{0}%
}\mathbf{L}_{b}\right)  =\left[  w\right]  \left(  \mathbf{L}_{a}%
w_{0}\mathbf{L}_{b}\right)  .
\]
The same argument (with the roles of $a$ and $b$ swapped) shows that%
\[
\left[  u\right]  \left(  \mathbf{L}_{b}\mathbf{L}_{a}\right)  =\left[
w\right]  \left(  \mathbf{L}_{b}w_{0}\mathbf{L}_{a}\right)  .
\]
The right-hand sides of these two equalities are equal, since%
\begin{align*}
\left[  w\right]  \left(  \mathbf{L}_{a}\underbrace{w_{0}\mathbf{L}_{b}%
}_{=b^{\ell\left(  w_{0}\right)  }\mathbf{L}_{1/b}}\right)   &  =b^{\ell
\left(  w_{0}\right)  }\underbrace{\left[  w\right]  \left(  \mathbf{L}%
_{a}\mathbf{L}_{1/b}\right)  }_{=\left[  w\right]  \left(  \mathbf{L}%
_{1/b}\mathbf{L}_{a}\right)  }=b^{\ell\left(  w_{0}\right)  }\left[  w\right]
\left(  \mathbf{L}_{1/b}\mathbf{L}_{a}\right) \\
&  =\left[  w\right]  \left(  \underbrace{b^{\ell\left(  w_{0}\right)
}\mathbf{L}_{1/b}}_{=\mathbf{L}_{b}w_{0}}\mathbf{L}_{a}\right)  =\left[
w\right]  \left(  \mathbf{L}_{b}w_{0}\mathbf{L}_{a}\right)  .
\end{align*}
Hence, so are the left-hand sides. That is, $\left[  u\right]  \left(
\mathbf{L}_{a}\mathbf{L}_{b}\right)  =\left[  u\right]  \left(  \mathbf{L}%
_{b}\mathbf{L}_{a}\right)  $. This shows that $u$ is $W$-length-symmetric.
This proves the \textquotedblleft only if\textquotedblright\ direction of
Proposition \ref{prop.Wls.ref}, and thus (as explained above) the whole proposition.
\end{proof}

\section{\label{sec.sep}Proofs: $W$-separability}

\subsection{Well-behavedness}

We shall next prove properties of $W$-separable elements. We begin with the
basic well-behavedness properties.

\begin{lemma}
\label{lem.sep1w0}Let $(W,S)$ be a Coxeter system with finite Coxeter group
$W$. Then, the elements $1$ and $w_{0}$ of $W$ are always $W$-separable.
\end{lemma}

\begin{proof}
For $1\in W$, this follows from Case 1 of Definition \ref{def.coxetersep}.

For $w_{0}\in W$, we argue as follows: If $S=\varnothing$, then this is clear
because $w_{0}=1$. Otherwise, pick any $s\in S$. The element $1\in
W_{S\setminus\left\{  s\right\}  }$ is $W_{S\setminus\left\{  s\right\}  }%
$-separable (by Case 1 of Definition \ref{def.coxetersep}). Hence, by Case 3
of Definition \ref{def.coxetersep}, the element $w_{0}1$ is $W$-separable. In
other words, $w_{0}$ is $W$-separable. This proves Lemma \ref{lem.sep1w0}.
\end{proof}

We extend the concept of an orthogonal decomposition to an arbitrary number of
subsets: An \emph{orthogonal decomposition} of a Coxeter system $\left(
W,S\right)  $ shall mean a tuple $\left(  J_{1},J_{2},\ldots,J_{k}\right)  $
of disjoint nonempty subsets of $S$ such that $S=J_{1}\cup J_{2}\cup\cdots\cup
J_{k}$ and such that all $i,j\in S$ that lie in different $J_{p}$'s satisfy
$m\left(  i,j\right)  =2$. We write such a decomposition as $S=\bigoplus
\limits_{p=1}^{k}J_{p}$. Again, any elements of $W_{J_{p}}$ for different
$p$'s commute with one another. If the only orthogonal decomposition of
$\left(  W,S\right)  $ is $\left(  S\right)  $ itself, then we say that the
Coxeter system $\left(  W,S\right)  $ is \emph{irreducible}. An orthogonal
decomposition $\left(  J_{1},J_{2},\ldots,J_{k}\right)  $ has maximum possible
$k$ if and only if all the Coxeter subsystems $\left(  W_{J_{p}},J_{p}\right)
$ are irreducible; in this case, the decomposition is called an
\emph{irreducible decomposition} of $\left(  W,S\right)  $. This decomposition
always exists (it can be obtained by picking the maximal irreducible subsets
$I$ of $S$), and is unique up to reordering. The subsets $J_{1},J_{2}%
,\ldots,J_{k}$ in this irreducible decomposition are called the
\emph{irreducible components} of $\left(  W,S\right)  $.

\begin{lemma}
\label{lem.w0-inv-irred}Let $\left(  W,S\right)  $ be an irreducible Coxeter
system with finite Coxeter group $W$. Let $p\in W$ be $W$-separable. Then,
$w_{0}p$ is also $W$-separable.
\end{lemma}

\begin{proof}
By Definition \ref{def.coxetersep}, the $W$-separability of $p$ means that we
are in one of the following four cases:

\textit{Case 1:} We have $p=1$\textbf{.}

\textit{Case 2:} The element $p$ is $W_{S\setminus\left\{  s\right\}  }%
$-separable for some $s\in S$.

\textit{Case 3:} We have $p=w_{0}u$, where $u\in W_{S\setminus\left\{
s\right\}  }$ is $W_{S\setminus\left\{  s\right\}  }$-separable for some $s\in
S$.

\textit{Case 4:} We have $p=uv$, where $S=I\oplus J$ is an orthogonal
decomposition of $\left(  W,S\right)  $, and where $u\in W_{I}$ and $v\in
W_{J}$ both are separable with respect to their respective parabolic subgroups.

In Case 1, we have $p=1$ and thus $w_{0}p=w_{0}$, which is $W$-separable by
Lemma \ref{lem.sep1w0}.

In Case 2, the element $w_{0}p$ is $W$-separable by Case 3 of Definition
\ref{def.coxetersep}, since $p\in W_{S\setminus\left\{  s\right\}  }$ is
$W_{S\setminus\left\{  s\right\}  }$-separable.

In Case 3, we have $p=w_{0}u$, so that $w_{0}p=w_{0}w_{0}u=u$ (since
$w_{0}w_{0}=1$), and this shows that $w_{0}p$ is $W_{S\setminus\left\{
s\right\}  }$-separable; therefore, $w_{0}p$ is $W$-separable as well by Case
2 of Definition \ref{def.coxetersep}.

Case 4 is impossible, since $W$ is irreducible and thus has no orthogonal
decomposition other than $\left(  S\right)  $.

So we are done in all cases.
\end{proof}

\begin{lemma}
\label{lem.orthogonalsep1}Let $(W,S)$ be a Coxeter system with finite Coxeter
group $W$ and with an orthogonal decomposition $S=\bigoplus\limits_{i=1}%
^{k}J_{i}$. For each $i\in\left\{  1,2,\ldots,k\right\}  $, let $w_{i}\in
W_{J_{i}}$ be a $W_{J_{i}}$-separable element. Then, the element $\prod
_{i=1}^{k}w_{i}\in W$ is $W$-separable.
\end{lemma}

\begin{proof}
We induct on $k$. Assume that $k>1$, since the other cases are trivial.
Consider the orthogonal decomposition $S=I\oplus J$ of $S$ where $I=J_{1}$ and
$J=\bigoplus_{i=2}^{k}J_{i}$ (this is nonempty since $k>1$). By induction, the
element $v=\prod_{i=2}^{k}w_{i}\in W_{J}$ is $W_{J}$-separable (since
$J=\bigoplus_{i=2}^{k}J_{i}$ is an orthogonal decomposition of $W_{J}$),
whereas $w_{1}\in W_{I}$ is $W_{I}$-separable (by assumption). Thus, by the
fourth case of Definition \ref{def.coxetersep}, the element $\prod_{i=1}%
^{k}w_{i}=w_{1}v$ is $W$-separable.
\end{proof}

\begin{lemma}
\label{lem.orthogonalsep} Let $(W,S)$ be a Coxeter system with finite Coxeter
group $W$ and with an irreducible decomposition $S=\bigoplus\limits_{i=1}%
^{k}J_{i}$. Let $w\in W$. Then $w$ is $W$-separable if and only if $w$ can be
written as $w=\prod_{i=1}^{k}w_{i}$ where $w_{i}\in W_{J_{i}}$ is $W_{J_{i}}%
$-separable for each $i$.
\end{lemma}

\begin{proof}
The backwards direction follows from Lemma \ref{lem.orthogonalsep1}.

For the forwards direction, we proceed by strong induction on $|S|$. We assume
that $w\in W$ is $W$-separable. By Definition \ref{def.coxetersep}, this can
be for four different reasons; we accordingly distinguish between four cases.

\textit{Case 1:} Suppose $w=1$. Then we can take $w_{i}=1$ for all $i$.

\textit{Case 2:} Suppose that $w\in W_{S\setminus\left\{  s\right\}  }$ is
$W_{S\setminus\{s\}}$-separable, where $s\in S$. Without loss of generality we
can assume that $s\in J_{1}$. Let $J_{1}\setminus\left\{  s\right\}
=\bigoplus_{\ell=1}^{r}I_{\ell}^{\prime}$ be an irreducible decomposition of
$W_{J_{1}\setminus\left\{  s\right\}  }$. Then, $S\setminus\left\{  s\right\}
=\left(  \bigoplus_{\ell=1}^{r}I_{\ell}^{\prime}\right)  \oplus\left(
\bigoplus_{i=2}^{k}J_{i}\right)  $ is an irreducible decomposition of
$W_{S\setminus\left\{  s\right\}  }$.

By induction, we thus have $w=(\prod_{\ell=1}^{r}v_{\ell})(\prod_{i=2}%
^{k}u_{i})$ where $v_{\ell}\in W_{I_{\ell}^{\prime}}$ is $W_{I_{\ell}^{\prime
}}$-separable for each $\ell$ and where $u_{i}\in W_{J_{i}}$ is $W_{J_{i}}%
$-separable for each $i$. Again by induction, $\prod_{\ell=1}^{r}v_{\ell}$ is
a $W_{J_{1}\setminus\{s\}}$-separable element of $W_{J_{1}\setminus\{s\}}$ and
thus a $W_{J_{1}}$-separable element of $W_{J_{1}}$ by Case 2 of Definition
\ref{def.coxetersep}. Consequently, if we let $w_{1}=\prod_{\ell=1}^{r}%
v_{\ell}$ and $w_{i}=u_{i}$ for all $i\in\left\{  2,3,\ldots,k\right\}  $,
then $w=\prod_{i=1}^{k}w_{i}$ where each $w_{i}$ is $W_{J_{i}}$-separable in
$W_{J_{i}}$. This again proves the claim.

\textit{Case 3:} Suppose that $w=w_{0}u$, where $u\in W_{S\setminus\left\{
s\right\}  }$ is $W_{S\setminus\left\{  s\right\}  }$-separable for some $s\in
S$.

For each $i\in\left\{  1,2,\ldots,k\right\}  $, let $w_{0,i}$ be the longest
element of the (finite) Coxeter subgroup $\left(  W_{J_{i}},J_{i}\right)  $.
Then, it is easy to see that $w_{0}=\prod_{i=1}^{k}w_{0,i}$ (see
(\ref{pf.prop.w0-inv.w0fac}) below).

But our above argument for Case 2 (applied to $u$ instead of $w$) shows that
$u$ can be written as $u=\prod_{i=1}^{k}u_{i}$ where each $u_{i}$ is
$W_{J_{i}}$-separable in $W_{J_{i}}$. Thus,%
\[
w=w_{0}u=\left(  \prod_{i=1}^{k}w_{0,i}\right)  \left(  \prod_{i=1}^{k}%
u_{i}\right)  =\prod_{i=1}^{k}\left(  w_{0,i}u_{i}\right)
\]
(since elements of distinct $W_{J_{i}}$'s commute). Moreover, for each
$i\in\left\{  1,2,\ldots,k\right\}  $, the element $w_{0,i}u_{i}\in W_{J_{i}}$
is $W_{J_{i}}$-separable by Lemma \ref{lem.w0-inv-irred} (applied to $\left(
W_{J_{i}},J_{i}\right)  $, $w_{0,i}$ and $u_{i}$ instead of $\left(
W,S\right)  $, $w_{0}$ and $p$), since $u_{i}$ is $W_{J_{i}}$-separable and
since $\left(  W_{J_{i}},J_{i}\right)  $ is irreducible. So we have written
$w$ as $w=\prod_{i=1}^{k}w_{i}$ where $w_{i}=w_{0,i}u_{i}\in W_{J_{i}}$ is
$W_{J_{i}}$-separable for each $i$. So our claim is proved here as well.

\textit{Case 4:} Suppose that $w=uv$ where $S=I\oplus J$ is an orthogonal
decomposition and $u\in W_{I}$ is $W_{I}$-separable and $v\in W_{J}$ is
$W_{J}$-separable. Let $I=\bigoplus_{\ell=1}^{r}I_{\ell}^{\prime}$ and
$J=\bigoplus_{\ell=1}^{t}J_{\ell}^{\prime}$ be irreducible decompositions of
$W_{I}$ and $W_{J}$, respectively. By induction, we can write $u=\prod
_{\ell=1}^{r}u_{\ell}$ where each $u_{\ell}\in W_{I_{\ell}^{\prime}}$ is
$W_{I_{\ell}^{\prime}}$-separable, and likewise write $v=\prod_{\ell=1}%
^{t}v_{\ell}$ where each $v_{\ell}$ is $W_{J_{\ell}^{\prime}}$-separable. This
proves our claim, because $w=uv=\left(  \prod_{\ell=1}^{r}u_{\ell}\right)
\left(  \prod_{\ell=1}^{t}v_{\ell}\right)  $ and since $S=\left(
\bigoplus_{\ell=1}^{r}I_{\ell}^{\prime}\right)  \oplus\left(  \bigoplus
_{\ell=1}^{t}J_{\ell}^{\prime}\right)  $ is an irreducible decomposition of
$W$ and therefore agrees with our given irreducible decomposition
$S=\bigoplus\limits_{i=1}^{k}J_{i}$ up to order.

The induction is thus complete in all four cases.
\end{proof}

We are now ready to prove Proposition \ref{prop.w0-inv}:

\begin{proof}
[Proof of Proposition \ref{prop.w0-inv}.]If $\left(  W,S\right)  $ is
irreducible, then this is Lemma \ref{lem.w0-inv-irred}. In the general case,
let $S=\bigoplus\limits_{i=1}^{k}J_{i}$ be an irreducible decomposition of
$\left(  W,S\right)  $. Then, the canonical map%
\begin{align*}
\prod_{i=1}^{k}W_{J_{i}}  &  \rightarrow W,\\
\left(  w_{i}\right)  _{i\in\left\{  1,2,\ldots,k\right\}  }  &  \mapsto
\prod_{i=1}^{k}w_{i}%
\end{align*}
is length-preserving (i.e., we have $\ell\left(  \prod_{i=1}^{k}w_{i}\right)
=\sum_{i=1}^{k}\ell\left(  w_{i}\right)  $ for all $\left(  w_{i}\right)
_{i\in\left\{  1,2,\ldots,k\right\}  }\in\prod_{i=1}^{k}W_{J_{i}}$). Thus, the
longest element $w_{0}$ of $W$ is given by
\begin{equation}
w_{0}=\prod_{i=1}^{k}w_{0,i}, \label{pf.prop.w0-inv.w0fac}%
\end{equation}
where each $w_{0,i}$ is the longest element of $W_{J_{i}}$.

Write $p$ as $p=\prod_{i=1}^{k}p_{i}$ where $p_{i}\in W_{J_{i}}$ for each $i$.
Then, the preceding two equalities yield
\[
w_{0}p=\left(  \prod_{i=1}^{k}w_{0,i}\right)  \left(  \prod_{i=1}^{k}%
p_{i}\right)  =\prod_{i=1}^{k}\left(  w_{0,i}p_{i}\right)
\]
(since elements of different $W_{J_{i}}$'s commute with each other).

Now, Lemma \ref{lem.orthogonalsep} shows that $p_{i}\in W_{J_{i}}$ is
$W_{J_{i}}$-separable for each $i$ (since $p\in W$ is $W$-separable, and since
$p=\prod_{i=1}^{k}p_{i}$). However, for each $i$, the Coxeter system $\left(
W_{J_{i}},J_{i}\right)  $ is irreducible, and thus Lemma
\ref{lem.w0-inv-irred} shows that $w_{0,i}p_{i}$ is $W_{J_{i}}$-separable
(since $p_{i}$ is $W_{J_{i}}$-separable). But Lemma \ref{lem.orthogonalsep}
(applied to $w_{0}p$ and $w_{0,i}p_{i}$ instead of $w$ and $w_{i}$) shows that
$w_{0}p$ is $W$-separable if and only if $w_{0,i}p_{i}\in W_{J_{i}}$ is
$W_{J_{i}}$-separable for each $i$. We thus conclude that $w_{0}p$ is
$W$-separable (since $w_{0,i}p_{i}\in W_{J_{i}}$ is $W_{J_{i}}$-separable for
each $i$).
\end{proof}

Our next goal is to show that separability is unchanged when we enlarge the
group, i.e., to prove Proposition \ref{prop.canonical}.

Broadly speaking, the proof will proceed by observing that it will suffice to
show the Proposition in the case where $(W,S)$ is irreducible. Once we have
done so, we can restrict to the intersection of $I$ within each component of
an irreducible decomposition. Before stating our first technical lemma, we
recall several auxiliary definitions associated to Coxeter groups. For more
depth we refer the reader to~\cite{BjoBre}. Firstly, let $T$ be the set of all
reflections of $\left(  W,S\right)  $ (that is, of all conjugates of elements
of $S$). Next, for $w\in W$, the \emph{support of} $w$, denoted by
$\operatorname*{supp}\left(  w\right)  $, is the set of simple reflections in
a reduced word for $w$. (This set does not depend on the reduced word, e.g.,
by Matsumoto's lemma.) If $\operatorname*{supp}\left(  w\right)  =S$, we say
that $w$ is \emph{fully supported}. Note that an element $w\in W$ and a subset
$I\subseteq S$ satisfy
\begin{equation}
w\in W_{I}\text{ if and only if }\operatorname*{supp}\left(  w\right)
\subseteq I. \label{eq.supp.WI}%
\end{equation}

\begin{lemma}
\label{lem.lower-int}Let $\left(  W,S\right)  $ be a Coxeter system. Let $I$
and $J$ be two subsets of $S$. Then, the set
\[
W_{I}W_{J}=\left\{  uv\ \mid\ u\in W_{I}\text{ and }v\in W_{J}\right\}
\]
is a lower ideal of the Bruhat order of $W$ (that is, if $p\in W_{I}W_{J}$ and
$q\leq p$, then $q\in W_{I}W_{J}$).
\end{lemma}

\begin{proof}
We must show that for any $u\in W_{I}$ and $v\in W_{J}$ and any $w\in W$
satisfying $w\leq uv$ (in the Bruhat order), we still have $w\in W_{I}W_{J}$.

Express $u$ and $v$ as $u=i_{1}i_{2}\cdots i_{k}$ and $v=j_{1}j_{2}\cdots
j_{\ell}$ with $i_{1},i_{2},\ldots,i_{k}\in I$ and $j_{1},j_{2},\ldots
,j_{\ell}\in J$. Thus, $uv=i_{1}i_{2}\cdots i_{k}j_{1}j_{2}\cdots j_{\ell}$.
By the deletion property, we can obtain a reduced word of $uv$ by removing
some factors from this product. That is, $uv$ has a reduced word of the form
\[
uv=i_{p_{1}}i_{p_{2}}\cdots i_{p_{x}}j_{q_{1}}j_{q_{2}}\cdots j_{q_{y}}%
\]
for some $p_{1},p_{2},\ldots,p_{x}\in\left\{  1,2,\ldots,k\right\}  $ and
$q_{1},q_{2},\ldots,q_{y}\in\left\{  1,2,\ldots,\ell\right\}  $. Since $w\leq
uv$, we can then obtain $w$ by removing some further factors from this reduced
word. The result is still a product of some elements of $I$ followed by some
elements of $J$; so it still belongs to $W_{I}W_{J}$. Thus, $w\in W_{I}W_{J}$, qed.
\end{proof}

\begin{lemma}
\label{lem.long_parabolic_darij}Let $\left(  W,S\right)  $ be an irreducible
Coxeter system. Let $I$ and $J$ be two proper subsets of $S$. Then, $W\neq
W_{I}W_{J}$.
\end{lemma}

\begin{proof}
Assume the contrary. Thus, $W=W_{I}W_{J}$.

Any element $w$ of $W$ can be written as a product of some $u\in W_{I}$ and
some $v\in W_{J}$ (since $w\in W=W_{I}W_{J}$). By the deletion property, we
can reduce this product to a reduced word by removing some of its factors;
hence, $w$ has a reduced word consisting entirely of elements of $I\cup J$.
Therefore,%
\begin{equation}
\operatorname*{supp}\left(  w\right)  \subseteq I\cup J.
\label{pf.lem.long_parabolic_darij.supp}%
\end{equation}

In particular, each $s\in S$ satisfies $\operatorname*{supp}\left(  s\right)
\subseteq I\cup J$ and therefore $s\in\operatorname*{supp}\left(  s\right)
\subseteq I\cup J$. Therefore, $I\cup J=S$.

Confronting this with $J\neq S$ (since $J$ is a proper subset of $S$), we find
$J\neq I\cup J$, so that $I\not \subseteq J$. Likewise, $J\not \subseteq I$.
Hence, there exist $i\in I\setminus J$ and $j\in J\setminus I$. Pick such $i$
and $j$.

Consider the \emph{Coxeter diagram} of $\left(  W,S\right)  $; this is the
undirected graph with vertex set $S$ and with an edge $ij$ for each
$2$-element subset $\left\{  i,j\right\}  $ of $S$ satisfying $m\left(
i,j\right)  >2$. The connected components of this diagram are precisely the
irreducible components of $\left(  W,S\right)  $. Since $\left(  W,S\right)  $
is irreducible, there must therefore be only one connected component. In other
words, the Coxeter diagram is connected. Thus, it has a path $\left(
v_{0},v_{1},\ldots,v_{r}\right)  $ from $j$ to $i$. That is, $v_{0}%
,v_{1},\ldots,v_{r}$ are distinct elements of $S$ with $v_{0}=j$ and $v_{r}%
=i$. Since $\left(  v_{0},v_{1},\ldots,v_{r}\right)  $ is a path, we have
$m\left(  v_{k},v_{k+1}\right)  >2$ for all $k\in\left\{  0,1,\ldots
,r-1\right\}  $.

Set $w:=v_{0}v_{1}\cdots v_{r}$. The word $v_{0}v_{1}\cdots v_{r}$ is reduced
and cannot be modified by any braid moves (since $m\left(  v_{k}%
,v_{k+1}\right)  >2$ for all $k\in\left\{  0,1,\ldots,r-1\right\}  $). Thus,
this word is the only reduced word of $w$ (since the Tits--Matsumoto theorem
shows that any two reduced words of $w$ can be transformed into each other
using braid moves).

However, $w$ can be written as a product of elements of $I$ multiplied with a
product of elements of $J$ (since $w\in W=W_{I}W_{J}$). By the deletion
condition, we can remove some factors from this product to obtain a reduced
word of $w$. This reduced word consists of some elements of $I$ followed by
some elements of $J$ (due to its origin), but it must be $v_{0}v_{1}\cdots
v_{r}$, since $v_{0}v_{1}\cdots v_{r}$ is the only reduced word of $w$.
However, this is impossible, because $v_{0}v_{1}\cdots v_{r}$ does
\textbf{not} consist of some elements of $I$ followed by some elements of $J$
(since $v_{0}=j\in J\setminus I$ and $v_{r}=i\in I\setminus J$). This
contradiction completes our proof.
\end{proof}

\begin{lemma}
\label{lem.long_parabolic} Let $(W,S)$ be an irreducible Coxeter system with
$W$ finite. Let $s\in S$ and $u\in W_{S\setminus\{s\}}$. Then $w_{0}u$ is
fully supported.
\end{lemma}

\begin{proof}
Assume the contrary. Then, some proper subset $I$ of $S$ satisfies
$\operatorname*{supp}\left(  w_{0}u\right)  =I$ and therefore $w_{0}u\in
W_{I}$. Therefore, $w_{0}\in W_{I}u^{-1}\subseteq W_{I}W_{S\setminus\left\{
s\right\}  }$ (because $u\in W_{S\setminus\left\{  s\right\}  }$).

Now, every element of $W$ is $\leq w_{0}$ and therefore also belongs to
$W_{I}W_{S\setminus\left\{  s\right\}  }$ (since Lemma \ref{lem.lower-int}
shows that $W_{I}W_{S\setminus\left\{  s\right\}  }$ is a lower ideal of the
Bruhat order of $W$). That is, $W=W_{I}W_{S\setminus\left\{  s\right\}  }$.
Moreover, $I$ and $S\setminus\left\{  s\right\}  $ are proper subsets of $S$.
Hence, Lemma \ref{lem.long_parabolic_darij} (for $J=S\setminus\left\{
s\right\}  $) shows that $W\neq W_{I}W_{S\setminus\left\{  s\right\}  }$. But
this contradicts $W=W_{I}W_{S\setminus\left\{  s\right\}  }$. This completes
the proof. \qedhere

\end{proof}

\begin{proof}
[Proof of Proposition~\ref{prop.canonical}]The backwards direction follows by
repeatedly applying Case 2 of Definition \ref{def.coxetersep} (strictly
speaking: by induction on $\left\vert S\right\vert $).

For the forwards direction, we again strongly induct on $\left\vert
S\right\vert $.

We will reduce to the case where $(W,S)$ is irreducible. So suppose that
$(W,S)$ is reducible with irreducible decomposition $S=\bigoplus_{i=1}%
^{k}J_{i}$ (so that $k>1$ and each $\left\vert J_{i}\right\vert <\left\vert
S\right\vert $). By Lemma~\ref{lem.orthogonalsep}, if $w$ is $W$-separable
then $w=\prod_{i=1}^{k}w_{i}$ where each $w_{i}\in W_{J_{i}}$ is $W_{J_{i}}%
$-separable. The orthogonality of the decomposition $S=\bigoplus_{i=1}%
^{k}J_{i}$ shows that elements of $W_{J_{i}}$ for distinct $i$'s commute, and
thus a reduced word for $w=\prod_{i=1}^{k}w_{i}$ can be obtained by
arbitrarily interleaving reduced words for all the $w_{i}$'s. Therefore,%
\[
\operatorname*{supp}\left(  w\right)  =\bigcup_{i=1}^{k}\operatorname*{supp}%
\left(  w_{i}\right)  .
\]
Hence, for each $i$, we have $\operatorname*{supp}\left(  w_{i}\right)
\subseteq\operatorname*{supp}\left(  w\right)  \subseteq I$ (by
(\ref{eq.supp.WI}), since $w\in W_{I}$). Combining this with
$\operatorname*{supp}\left(  w_{i}\right)  \subseteq J_{i}$ (since $w_{i}\in
W_{J_{i}}$), we obtain $\operatorname*{supp}\left(  w_{i}\right)  \subseteq
I\cap J_{i}$ and therefore $w_{i}\in W_{I\cap J_{i}}$ for each $i$. Hence, by
induction on $|S|$, we must have that $w_{i}$ is $W_{I\cap J_{i}}$-separable
for each $i$ (since $w_{i}$ is $W_{J_{i}}$-separable and $\left\vert
J_{i}\right\vert <\left\vert S\right\vert $). Since $I=\bigoplus_{i=1}%
^{k}\left(  I\cap J_{i}\right)  $ is an orthogonal decomposition for $W_{I}$
(after omitting empty intersections), we thus conclude that $w$ is $W_{I}%
$-separable by Lemma \ref{lem.orthogonalsep1}.

We now turn our attention to the case where $(W,S)$ is irreducible. We assume
that $I\neq S$, since otherwise the claim is obvious. We have $w\in W_{I}$,
thus $\operatorname*{supp}\left(  w\right)  \subseteq I$ by (\ref{eq.supp.WI}%
). Hence, $w$ cannot be fully supported (since $I\neq S$).

By Definition \ref{def.coxetersep}, since $w\in W$ is separable, we must be in
one of the four cases listed in the definition. But we cannot be in Case 4,
since the irreducibility of $W$ precludes an orthogonal decomposition. Nor can
we be in Case 3, since (again by irreducibility of $W$) Lemma
\ref{lem.long_parabolic} would force $w$ to be fully supported in this case.
So we must be in Case 1 or Case 2. In Case 1, we have $w=1$ and thus the claim
holds trivially. Assume that we are in Case 2. Thus, $w$ is a $W_{S\setminus
\{s\}}$-separable element of $W_{S\setminus\{s\}}$ for some $s\in S$. So
(\ref{eq.supp.WI}) yields $\operatorname*{supp}\left(  w\right)  \subseteq
S\setminus\left\{  s\right\}  $. Combined with $\operatorname*{supp}\left(
w\right)  \subseteq I$, this yields $\operatorname*{supp}\left(  w\right)
\subseteq\left(  S\setminus\left\{  s\right\}  \right)  \cap I=I\setminus
\left\{  s\right\}  $. By induction on $\left\vert S\right\vert $, we now
conclude that $w$ is $W_{I\setminus\left\{  s\right\}  }$-separable (since $w$
is $W_{S\setminus\left\{  s\right\}  }$-separable). Hence, it follows that $w$
is $W_{I}$-separable as well (tautologically if $s\notin I$; otherwise by Case
2 in Definition \ref{def.coxetersep}). This proves our claim, thus completing
the induction.
\end{proof}

With Proposition \ref{prop.canonical} proved, we now have the right to drop
the \textquotedblleft$W$\textquotedblright\ in the word \textquotedblleft%
$W$-separable\textquotedblright: The $W$-separability of an element $w\in W$
still depends on the Coxeter system $\left(  W,S\right)  $, but at least does
not change based on moving to a (standard parabolic) Coxeter subsystem.

\subsection{Separable elements in type A}

\begin{proof}
[Proof of Theorem \ref{thm.separable_equivalence}.]In this proof, the word
\textquotedblleft separable\textquotedblright\ without a prefix shall always
mean \textquotedblleft separable as a permutation\textquotedblright, whereas
$W$-separability for a Coxeter group $W$ will always carry the
\textquotedblleft$W$-\textquotedblright\ prefix.

We recall the four operations $\operatorname*{rev}$, $\operatorname*{com}$,
$\oplus$ and $\ominus$ on permutations defined as follows:

\begin{enumerate}
\item The operation $\operatorname*{rev}$ (called \textquotedblleft
reversal\textquotedblright) sends any permutation $w\in S_{k}$ to the
permutation $\operatorname*{rev}w:=ww_{0}\in S_{k}$, where $w_{0}$ is the
longest element of $S_{k}$. (So $\left(  \operatorname*{rev}w\right)  \left(
i\right)  =w\left(  k+1-i\right)  $ for each $i\in\left[  k\right]  $. In
one-line notation, $\operatorname*{rev}w$ is simply obtained by reversing $w$.)

\item The operation $\operatorname*{com}$ (called \textquotedblleft
complementation\textquotedblright) sends any permutation $w\in S_{k}$ to the
permutation $\operatorname*{com}w:=w_{0}w\in S_{k}$, where $w_{0}$ is the
longest element of $S_{k}$. (So $\left(  \operatorname*{com}w\right)  \left(
i\right)  =k+1-w\left(  i\right)  $ for each $i\in\left[  k\right]  $.)

\item The operation $\oplus$ (called \textquotedblleft direct
sum\textquotedblright) sends two permutations $u\in S_{k}$ and $v\in S_{\ell}$
to the permutation $u\oplus v\in S_{k+\ell}$ that acts as $u$ on the numbers
$1,2,\ldots,k$ and as (a shifted copy of) $v$ on the numbers $k+1,k+2,\ldots
,k+\ell$.

\item The operation $\ominus$ (called \textquotedblleft skew
sum\textquotedblright) sends two permutations $u\in S_{k}$ and $v\in S_{\ell}$
to the permutation $u\ominus v:=\operatorname*{rev}\left(  \left(
\operatorname*{rev}u\right)  \oplus\left(  \operatorname*{rev}v\right)
\right)  \in S_{k+\ell}$. (Somewhat confusingly, $u\ominus v$ is denoted
$v\ominus u$ in \cite{Kitaev11}; but this changes nothing about separability.)
\end{enumerate}

Among the many characterizations of separable permutations given in
\cite[\textquotedblleft From a pruned tree to a polynomial interchange and a
separable permutation\textquotedblright, Theorem]{Ghys16}, we will need the
fourth one (see also \cite[Definition 2.2.34]{Kitaev11}): A permutation $w\in
S_{n}$ is separable if and only if it can be obtained from several copies of
the trivial permutation on one object (i.e., the permutation
$\operatorname*{id}\in S_{1}$) by successive $\oplus$ and $\ominus$
operations. In other words, a permutation $\sigma\in S_{n+1}$ is separable if
and only if $\sigma$ is either the identity permutation or can be factored as
$\sigma=\sigma_{1}\oplus\sigma_{2}$ or $\sigma=\sigma_{1}\ominus\sigma_{2}$
where both $\sigma_{1}\in S_{k}$ and $\sigma_{2}\in S_{n+1-k}$ are separable
(and where $1\leq k\leq n$).

We also note that the set of all separable permutations is fixed under both
operations $\operatorname*{rev}$ and $\operatorname*{com}$. This is easy to
see using the pattern-avoidance characterization of separable permutations,
but also follows from the recursive definition just recalled, after observing
that any permutations $u\in S_{k}$ and $v\in S_{\ell}$ satisfy
\begin{align*}
\operatorname*{rev}\left(  u\ominus v\right)   &  =\left(  \operatorname*{rev}%
u\right)  \oplus\left(  \operatorname*{rev}v\right)
,\ \ \ \ \ \ \ \ \ \ \operatorname*{rev}\left(  u\oplus v\right)  =\left(
\operatorname*{rev}u\right)  \ominus\left(  \operatorname*{rev}v\right)  ,\\
\operatorname*{com}\left(  u\ominus v\right)   &  =\left(  \operatorname*{com}%
v\right)  \oplus\left(  \operatorname*{com}u\right)
,\ \ \ \ \ \ \ \ \ \ \operatorname*{com}\left(  u\oplus v\right)  =\left(
\operatorname*{com}v\right)  \ominus\left(  \operatorname*{com}u\right)  .
\end{align*}
The set of all $A_{n}$-separable elements of $A_{n}$ is also preserved under
$\operatorname*{com}$ by Proposition \ref{prop.w0-inv} (or even by Lemma
\ref{lem.w0-inv-irred}, since the Coxeter system $A_{n}$ is irreducible).

Write $W$ for the Coxeter group $A_{n}=S_{n+1}$.

Now we have to prove that an element $u\in A_{n}$ is $A_{n}$-separable if and
only if $u$ is separable.

\begin{enumerate}
\item[$\Longrightarrow:$] For the forward direction, we proceed by strong
induction on $n$. Note that the base case of $A_{0}$ is trivial, since
$A_{0}=\left\{  1\right\}  $. Now assume that our claim holds for all $A_{i}$
with $0\leq i<n$, and consider an $A_{n}$-separable element $w\in A_{n}$.
Since $A_{n}$ is irreducible, there is no nontrivial orthogonal decomposition;
so we need only verify for each of the remaining three possibilities of
Definition~\ref{def.coxetersep} that $w$ is separable as a permutation in
$S_{n+1}$.

In the first case where $w=1$, this follows from the fact that the identity
permutation $1=\operatorname*{id}$ is separable (as it is the direct sum of
$n+1$ copies of the one-element permutation).

Consider the second case, where $w$ is a $W_{S\setminus\{s\}}$-separable
element of a parabolic subgroup $W_{S\setminus\{s\}}$. Specifically, let $w$
be a $W_{S\setminus\{s_{k}\}}$-separable element of $W_{S\setminus\{s_{k}\}}$
for $k\in\lbrack n]$.

Let $I=\{s_{1},\dots,s_{k-1}\}$ and let $J=\{s_{k+1},\dots,s_{n}\}$. Then,
$S\setminus\{s_{k}\}=I\oplus J$ is an irreducible decomposition\footnote{This
is a bit false for $k=1$ and for $k=n$; in these cases, the argument is even
simpler and left to the reader.} of $S\setminus\{s_{k}\}$. Hence,
Lemma~\ref{lem.orthogonalsep} shows that $w=uv$ where $u\in W_{I}$ is $W_{I}%
$-separable and $v\in W_{J}$ is $W_{J}$-separable. The parabolic subgroups
$W_{I}$ and $W_{J}$ are canonically isomorphic (as Coxeter groups) to the
symmetric groups $S_{k}$ and $S_{n+1-k}$, respectively; the elements $u$ and
$v$ of the former correspond to two permutations $\sigma\in S_{k}$ and
$\tau\in S_{n+1-k}$ in the latter. The equality $w=uv$ then becomes
$w=\sigma\oplus\tau$ (as a permutation). Moreover, the elements $\sigma$ and
$\tau$ are (respectively) $A_{k-1}$-separable and $A_{n-k}$-separable (since
$u\in W_{I}$ is $W_{I}$-separable and $v\in W_{J}$ is $W_{J}$-separable), and
thus -- by our induction hypothesis -- are separable permutations. Hence,
their direct sum $\sigma\oplus\tau$ is also separable by the above definition
of separability. That is, $w$ is separable as a permutation.

Finally, consider the third case, in which $w=w_{0}u$ where $u$ is a
$W_{S\setminus\{s\}}$-separable element of a parabolic subgroup $W_{S\setminus
\{s\}}$. By our argument in the previous case, $u$ is separable as a
permutation. Hence, $\operatorname*{com}u$ is separable as a permutation as
well, since the set of all separable permutations is preserved under
$\operatorname*{com}$. That is, $w$ is a separable permutation, since
$w=w_{0}u=\operatorname*{com}u$.

Thus, the forward direction of Theorem \ref{thm.separable_equivalence} is proved.

\item[$\Longleftarrow:$] For the backward direction, we again proceed by
strong induction. The base case of $1\in S_{1}$ holds trivially for the same
reasons as before.

Let us now fix $n\geq1$, and assume that every separable permutation $\pi\in
S_{i}$ is $A_{i-1}$-separable for each $i\leq n$. Next let $\sigma\in S_{n+1}$
be a separable permutation. We must show that $\sigma$ is $A_{n}$-separable.

Since $\sigma$ is a separable permutation, we have $\sigma=\sigma_{1}%
\oplus\sigma_{2}$ or $\sigma=\sigma_{1}\ominus\sigma_{2}$, where $\sigma
_{1}\in S_{k}$ and $\sigma_{2}\in S_{n-k+1}$ are both separable permutations.
We WLOG assume that $\sigma=\sigma_{1}\oplus\sigma_{2}$, since the case of
$\sigma=\sigma_{1}\ominus\sigma_{2}$ can be reduced to the former case using
$\operatorname*{com}\left(  u\ominus v\right)  =\left(  \operatorname*{com}%
v\right)  \oplus\left(  \operatorname*{com}u\right)  $ (and using the fact
that the set of all separable permutations and the set of all $A_{n}%
$-separable elements of $A_{n}$ are preserved under $\operatorname*{com}$). In
the language of Coxeter groups, this means that $\sigma=uv$ for some $u\in
W_{I}$ and some $v\in W_{J}$, where $I=\left\{  s_{1},s_{2},\ldots
,s_{k-1}\right\}  $ and $J=\left\{  s_{k+1},s_{k+2},\ldots,s_{n}\right\}  $
(namely, $u$ is the element of $W_{I}$ corresponding to $\sigma_{1}$ under the
canonical isomorphism $S_{k}=A_{k-1}\cong W_{I}$, whereas $v$ is the element
of $W_{J}$ corresponding to $\sigma_{2}$ under the canonical isomorphism
$S_{n-k+1}=A_{n-k}\cong W_{J}$). Moreover, these elements $u$ and $v$
correspond to the separable permutations $\sigma_{1}$ and $\sigma_{2}$ in
$S_{k}$ and $S_{n-k+1}$. By the induction hypothesis, the latter separable
permutations $\sigma_{1}$ and $\sigma_{2}$ are (respectively) $A_{k-1}%
$-separable and $A_{n-k}$-separable (since $k$ and $n-k+1$ are $\leq n$).
Hence, the corresponding elements $u$ and $v$ are (respectively) $W_{I}%
$-separable and $W_{J}$-separable (because of the Coxeter group isomorphisms
$A_{k-1}\cong W_{I}$ and $A_{n-k}\cong W_{J}$). Their product $uv$ is
therefore $W_{S\setminus\left\{  s_{k}\right\}  }$-separable, since
$S\setminus\left\{  s_{k}\right\}  =I\oplus J$ is an orthogonal decomposition
of $\left(  W_{S\setminus\left\{  s_{k}\right\}  },S\setminus\left\{
s_{k}\right\}  \right)  $. Hence, $uv$ is also $W$-separable (by the
definition of $W$-separability). In other words, $uv$ is $A_{n}$-separable. In
other words, $\sigma$ is $A_{n}$-separable (since $\sigma=uv$). This completes
the proof of the backward direction. \qedhere

\end{enumerate}
\end{proof}

\subsection{Length-symmetry of separable elements}

\begin{proof}
[Proof of Theorem \ref{thm.sep}.]We prove this by strong induction on
$\left\vert S\right\vert $, using the recursive definition of the separability
of $w$.

\begin{itemize}
\item If $w=1$, then $w$ is $W$-length-symmetric by
Proposition~\ref{prop.Wls.inv}, as the identity element is an involution.

\item If $w\in W_{S\setminus\left\{  s\right\}  }$ is $W_{S\setminus\left\{
s\right\}  }$-separable for some $s\in S$, then we know by the induction
hypothesis that $w$ is $W_{S\setminus\left\{  s\right\}  }$-length-symmetric.
Hence, by Proposition \ref{prop.Wls.WI}, we conclude that $w$ is $W$-length-symmetric.

\item If $w=w_{0}u$, where $u$ is a separable element of a parabolic subgroup
$W_{S\setminus\left\{  s\right\}  }$ for some simple reflection $s\in S$, then
the case we have just done shows that $u$ is $W$-length-symmetric. Hence,
Proposition \ref{prop.Wls.ref} shows that $w_{0}u$ is $W$-length-symmetric. In
other words, $w$ is $W$-length-symmetric.

\item If $w=uv$, where $S=I\oplus J$ is an orthogonal decomposition of $W$,
and where $u\in W_{I}$ and $v\in W_{J}$ both are separable (with respect to
$W_{I}$ and $W_{J}$), then the induction hypothesis shows that $u$ is $W_{I}%
$-length-symmetric and that $v$ is $W_{J}$-length-symmetric. Thus, Proposition
\ref{prop.Wls.I+J} shows that $uv$ is $W$-length-symmetric. In other words,
$w$ is $W$-length-symmetric.
\end{itemize}

So the claim that $w$ is $W$-length-symmetric is proved in all possible cases.
\end{proof}

\subsection{\label{subsec.sep.gaga}Equivalence with the Gaetz--Gao definition}

Previously in~\cite[Definition 3.1]{GaeGao20a}, \cite[Definition 1]{GaeGao20b}
and \cite[Definition 5.1.1]{Gaetz21}, Gaetz and Gao introduced a
generalization of separable permutations to Weyl groups. Formally, they
defined separable elements of a Weyl group as follows. For a discussion of
root systems and their associated terminology, we refer the reader
to~\cite[Chapters 1, 2, 5]{Humphreys} and \cite[\S 9, \S 10]{Humphreys72}. Let
$\Phi$ be a root system with set of simple roots $\Delta$ and positive roots
$\Phi^{+}$. Let $W(\Phi)$ denote the associated Weyl group. The
\emph{inversion set} of an element $w\in W(\Phi)$ is defined to be the set%
\[
I_{\Phi}(w)=\left\{  \alpha\in\Phi^{+}\mid w\alpha\in\Phi^{-}\right\}  .
\]
If $w\in W\left(  \Phi\right)  $ is any element and $\Phi^{\prime}$ is a root
subsystem of $\Phi$ generated by a subset $\Delta^{\prime}$ of $\Delta$, then
the \emph{restriction} $w|_{\Phi^{\prime}}$ of $w$ to $\Phi^{\prime}$ is
defined as the unique element $w^{\prime}$ of $W\left(  \Phi^{\prime}\right)
$ satisfying $I_{\Phi^{\prime}}\left(  w^{\prime}\right)  =I_{\Phi}\left(
w\right)  \cap\left(  \Phi^{\prime}\right)  ^{+}$. The set $\Phi$ is partially
ordered by the \emph{root order}, in which $\beta\geq\alpha$ if and only if
$\beta-\alpha$ is a sum of positive roots.

\begin{definition}
[cf.~\cite{GaeGao20a}]\label{def.ggseparable} Let $w\in W(\Phi)$. Then $w$ is
\textit{separable} if one of the following holds:

\begin{itemize}
\item $\Phi$ is of type $A_{1}$;

\item $\Phi=\bigoplus\limits_{i}\Phi_{i}$ is reducible and $w|_{\Phi_{i}}$ is
separable for each $i$;

\item $\Phi$ is irreducible and there exists a simple root $\alpha_{i}%
\in\Delta$ such that $w|_{\Phi^{\prime}}\in W(\Phi^{\prime})$ is separable,
where $\Phi^{\prime}$ is the root system generated by $\Delta^{\prime}%
=\Delta\setminus\{\alpha_{i}\}$, and such that
\begin{align*}
\text{either\ }  &  \{\beta\in\Phi^{+}:\beta\geq\alpha_{i}\}\subseteq I_{\Phi
}(w),\\
\text{or\ }  &  \{\beta\in\Phi^{+}:\beta\geq\alpha_{i}\}\cap I_{\Phi
}(w)=\varnothing.
\end{align*}
Such a root $\alpha_{i}$ is called a \emph{pivot} for $w$.
\end{itemize}
\end{definition}

As we shall immediately see, these two definitions are in fact equivalent when
$W$ is a Weyl group.

\begin{theorem}
\label{thm.gg} Let $(W,S)$ be a finite Weyl group. Then an element $w\in W$ is
$W$-separable as in Definition~\ref{def.coxetersep} if and only if $w$ is
separable as in Definition~\ref{def.ggseparable}.
\end{theorem}

\begin{proof}
$\Longrightarrow:$ Suppose that $w\in W$ is separable as in
Definition~\ref{def.coxetersep}. We show that $w$ will satisfy at least one of
the conditions of Definition~\ref{def.ggseparable}.

We use strong induction on the rank of $\left(  W,S\right)  $ (that is, the
size of $S$).

Suppose first that $W$ is reducible: $S=\bigoplus_{i}S_{i}$, with
corresponding decomposition $\Phi=\bigoplus\limits_{i}\Phi_{i}$ of root
systems. By Lemma~\ref{lem.orthogonalsep}, we have that $w$ is a product of
separable elements\footnote{By the induction hypothesis, we need not
distinguish between the two definitions of separability for any Coxeter system
of rank $<\left\vert S\right\vert $.} of the $W_{S_{i}}$. These separable
elements are precisely the $w|_{\Phi_{i}}$, so we conclude that $w$ is
separable in the sense of Definition~\ref{def.ggseparable}. (This includes the
case $S=\varnothing$, in which case $w=1$ is clearly separable.)

Assume now that $W$ is irreducible and $w$ is separable with $w\in
W_{S\setminus\{s\}}$ separable for some simple generator $s$. Since $w\in
W_{S\setminus\{s\}}$, no positive root $\beta\geq\alpha_{s}$ (where
$\alpha_{s}$ is the root corresponding to the reflection $s$) is sent to a
negative root under the action of $w$; that is, $I_{\Phi}(w)\cap\{\beta\in
\Phi^{+}:\beta\geq\alpha_{s}\}=\varnothing$. Hence, if we let $\Phi^{\prime}$
denote the root system generated by $\Delta^{\prime}=\Delta\setminus
\{\alpha_{s}\}$, then $w|_{\Phi^{\prime}}=w$ (since these two elements have
the same inversion set). By induction, $w|_{\Phi^{\prime}}\in W(\Phi^{\prime
})=W_{S\setminus\{s\}}$ is separable. So $w$ is separable in the sense of
Definition~\ref{def.ggseparable}.

For the case $w=w_{0}u$ with $u\in W_{S\setminus\{s\}}$ separable, we have
$I_{\Phi}(w_{0}u)=\Phi^{+}\setminus I_{\Phi}(u)$ since $I_{\Phi}(w_{0}%
)=\Phi^{+}$. As $I_{\Phi}(u)\cap\{\beta\in\Phi^{+}:\beta\geq\alpha
_{s}\}=\varnothing$ by the previous case, $\alpha_{s}$ is a pivot for $w$.
Moreover, letting $w_{0}^{\prime}$ denote the longest element of
$W_{S\setminus\{s\}}$, we have $I_{\Phi}(w)\cap(\Phi^{\prime})^{+}%
=(\Phi^{\prime})^{+}\setminus I_{\Phi^{\prime}}(u)=I_{\Phi^{\prime}}%
(w_{0}^{\prime}u)$, so $w|_{\Phi^{\prime}}=w_{0}^{\prime}u$. By induction $u$
is separable in $W(\Phi^{\prime})$, hence so is $w_{0}^{\prime}u$ by
\cite[Proposition 3.5]{GaeGao20a} applied to $\Phi^{\prime}$. Thus $w$ is
separable in the sense of Definition~\ref{def.ggseparable} (since
$w|_{\Phi^{\prime}}=w_{0}^{\prime}u$).

$\Longleftarrow:$ Now suppose that $w$ is separable as in
Definition~\ref{def.ggseparable}. Note that any element of type $A_{1}$
includes $1$ so our base cases are the same. If $(W,S)$ is reducible, the
second condition of Definition~\ref{def.ggseparable} is precisely
Lemma~\ref{lem.orthogonalsep}, since $w=\prod_{i}\left(  w|_{\Phi_{i}}\right)
$. Assume now that $(W,S)$ is irreducible and $w$ has a pivot $\alpha_{i}$
such that $w|_{\Phi^{\prime}}\in W(\Phi^{\prime})$ is separable for
$\Delta^{\prime}=\Delta\setminus\{\alpha_{i}\}$. By the definition of a pivot,
we have%
\begin{align*}
\text{either\ }  &  \{\beta\in\Phi^{+}:\beta\geq\alpha_{i}\}\subseteq I_{\Phi
}(w),\\
\text{or\ }  &  \{\beta\in\Phi^{+}:\beta\geq\alpha_{i}\}\cap I_{\Phi
}(w)=\varnothing.
\end{align*}
In the latter case, $w=w|_{\Phi^{\prime}}$ (since both sides have the same
inversion set), so that $w$ is a separable element of $W\left(  \Phi^{\prime
}\right)  =W_{S\setminus\{s_{\alpha_{i}}\}}$ where $s_{\alpha_{i}}$ is the
simple reflection corresponding to $\alpha_{i}$, and therefore $w$ is also
$W$-separable in the sense of Definition~\ref{def.coxetersep}. In the former
case, let $u=w_{0}w$. By \cite[Proposition 3.5]{GaeGao20a}, the element $u$ is
also separable, and $u|_{\Phi^{\prime}}$ is separable by \cite[Proposition
5.2]{GaeGao20a}. Further, since $\{\beta\in\Phi^{+}:\beta\geq\alpha_{i}\}\cap
I_{\Phi}(u)=\varnothing$, then, just as in the previous case, $u$ is a
separable element of $W_{S\setminus\{s_{\alpha_{i}}\}}$ where $s_{\alpha_{i}}$
is the simple reflection corresponding to $\alpha_{i}$. Thus we can conclude
that $w=w_{0}u$ where $u$ is a separable element of $W_{S\setminus\{s\}}$ for
some $s$; therefore, $w$ is $W$-separable (in the sense of
Definition~\ref{def.coxetersep}) as well, proving the claim.
\end{proof}

\subsection{\label{subsec.sep.count}Counting separable elements}

We now turn our attention to the enumeration of separable elements in finite
Coxeter systems. In the type-A case, the answer is well-known (see, e.g.,
\cite[\textquotedblleft Hipparchus and Schroeder\textquotedblright]{Ghys16},
\cite[Theorem 4.3]{West95}, \cite[Exercise 6.39]{Stanley-EC2}): The number of
separable permutations in $S_{n}$ is the $\left(  n-1\right)  $-st \emph{large
Schr\"{o}der number}.

Next we handle the easy case of the dihedral groups:

\begin{proposition}
\label{prop.dihedral_sep_enum} Let $(W,S)$ be a Coxeter system with finite
Coxeter group $W$ such that $S=\left\{  s_{1},s_{2}\right\}  $ and $s_{1}%
s_{2}\neq s_{2}s_{1}$; that is, let $W$ be a finite dihedral group with
connected Coxeter diagram. Then the separable elements of $W$ are
$1,s_{1},s_{2},w_{0}s_{1},w_{0}s_{2},w_{0}$.
\end{proposition}

\begin{proof}
We first verify that the listed elements are indeed separable. First note that
$S$ has no orthogonal decomposition (since the Coxeter diagram is connected),
and thus all separable elements $w$ are either the identity, $w\in
W_{\{s_{i}\}}$ where $w$ is separable, or $w_{0}w\in W_{\{s_{i}\}}$ where
again $w$ is separable. But every element in a Coxeter system with a single
generator is separable. Consequently we must then have that $1,s_{1}%
,s_{2},w_{0}s_{1},w_{0}s_{2},w_{0}$ are all separable. What remains to be seen
is that there are no other separable elements. If $w$ is separable and
$w\notin\{1,w_{0}\}$ then $w$ or $w_{0}w$ is a nontrivial element of a
parabolic subgroup. But the only such elements are $s_{1},s_{2}$, so we have
found all separable elements as claimed.
\end{proof}

We let $\sep\left(  W,S\right)  $ or $\sep\left(  W\right)  $ denote the
number of all separable elements of a Coxeter system $(W,S)$. This number is
multiplicative with respect to orthogonal decompositions:

\begin{lemma}
\label{lem.enumorthog} Let $(W,S)$ be a Coxeter system with finite Coxeter
group $W$ and suppose that $(I,J)$ is an orthogonal decomposition of $(W,S)$.
Then $\sep(W,S)=\sep(W_{I},I)\cdot\sep(W_{J},J)$.
\end{lemma}

\begin{proof}
Let $\Sep(W,S)$ denote the set of all separable elements of a Coxeter system
$(W,S)$. Lemma~\ref{lem.orthogonalsep} shows that there is a bijection
\begin{align*}
\Sep(W_{I},I)\times\Sep(W_{J},J)  &  \rightarrow\Sep(W,S),\\
\left(  u,v\right)   &  \mapsto uv.
\end{align*}
Hence, by the bijection principle, $\left\vert \Sep(W_{I},I)\times
\Sep(W_{J},J)\right\vert =\left\vert \Sep(W,S)\right\vert $. That is,
$\sep(W,S)=\sep(W_{I},I)\cdot\sep(W_{J},J)$.
\end{proof}

\begin{lemma}
\label{lem.enumsep}Let $(W,S)$ be an irreducible Coxeter system with finite
Coxeter group $W$. Then,%
\[
\sep(W)=2\sum_{I\subsetneq S}(-1)^{|S|-|I|-1}\sep(W_{I})
\]

\end{lemma}

\begin{proof}
If $w\in W$ is $W$-separable, then $w$ must satisfy one of the first three
cases in Definition \ref{def.coxetersep} (since the irreducibility of $\left(
W,S\right)  $ renders Case 4 impossible). In the first two cases, we have
$w\in W_{S\setminus\left\{  s\right\}  }$ for some $s\in S$ (in Case 1, any
$s$ will do), which shows that $w$ is not fully supported, whereas Lemma
\ref{lem.long_parabolic} shows that $w_{0}w$ is fully supported (since $w\in
W_{S\setminus\left\{  s\right\}  }$). On the other hand, in Case 3, we have
$w=w_{0}u$ for some $s\in S$ and $u\in W_{S\setminus\left\{  s\right\}  }$,
which shows that $w_{0}w=\underbrace{w_{0}w_{0}}_{=1}u=u$ is not fully
supported (since $u\in W_{S\setminus\left\{  s\right\}  }$), whereas
$w=w_{0}u$ is fully supported (by Lemma \ref{lem.long_parabolic}). Thus, in
each case, we have shown that \textbf{exactly one} of $w$ and $w_{0}w$ is
fully supported. Thus, the $W$-separable elements $w\in W$ can be classified
into two (disjoint) types: the ones that are fully supported, and the ones for
which $w_{0}w$ is fully supported. Left multiplication by $w_{0}$ is a
bijection between these two classes.

Hence, it suffices to show that the number of separable elements of $(W,S)$
which are not fully supported is $\sum_{I\subsetneq S}(-1)^{|S|-|I|-1}%
\sep(W_{I})$. This is because a separable element of $W$ which is not fully
supported on $S$ is precisely a separable element of $W$ that is contained
within a standard parabolic subgroup $W_{I}$ of $W$ with $I\subsetneq S$. Thus
the claimed enumeration will follow by inclusion-exclusion.
\end{proof}

\begin{proposition}
The number of separable elements of $B_{n}$ is enumerated by the large
Schr\"{o}der numbers. That is, $\sep(B_{n})=\sep(A_{n})$.
\end{proposition}

\begin{proof}
We proceed by strong induction on $n$. For a base case we consider $n=2$, as
this is the minimal $n$ in which $B_{n}$ and $A_{n}$ are distinct. Then,
$\sep(B_{2})=\sep(A_{2})$ follows from
Proposition~\ref{prop.dihedral_sep_enum}. So now let $n\geq3$, and assume that
for each $k<n$ we have that $\sep(B_{k})=\sep(A_{k})$. Observe that every
proper parabolic subgroup of $B_{n}$ has irreducible components of type
$A_{i}$ for some $i$'s as well as at most one component of type $B_{j}$ for
some $j$. Importantly, for any such subgroup, there is a corresponding
subgroup of $A_{n}$ but where the $B_{j}$-factor is instead $A_{j}$; this
latter subgroup has the same Coxeter diagram as the former except for possibly
an edge label being $3$ instead of $4$. By the induction hypothesis, the
latter subgroup of $A_{n}$ will have the same number of separable elements as
the former subgroup of $B_{n}$ (since Lemma \ref{lem.enumorthog} shows that
both numbers equal the products of the corresponding numbers for the
irreducible components, which are $A_{k}$s or $B_{k}$s). Hence, by
Lemma~\ref{lem.enumsep}, we conclude that the number of separable elements of
$B_{n}$ will satisfy the same recursion as $A_{n}$ with the same initial conditions.
\end{proof}

\section{Further questions}

As Theorem \ref{thm.sep} shows, the $W$-separable elements are $W$%
-length-symmetric, but this is not an \textquotedblleft if and only
if\textquotedblright\ in general. One can easily construct larger classes of
$W$-length-symmetric elements by adding further production rules to the
recursive definition of $W$-separability; for example, all involutions could
be included in the set at its onset. The resulting classes of elements of $W$
are still likely to not exhaust all $W$-length-symmetric elements, but they
might have interesting properties of their own; in particular, one can try to
enumerate them for $W=S_{n}$.

On the other hand, as we mentioned in Remark \ref{rmk.hecke}, there is a
stronger notion of Hecke-$W$-length-symmetry that warrants further study and enumeration.

Our proofs of the $W$-length-symmetry results have a rather transparent
computational structure. One might wonder whether they could lead to
recursions for $\left[  w\right]  \left(  \mathbf{L}_{a}\mathbf{L}_{b}\right)
$ that make the $\left(  a,b\right)  $-symmetry visible when $w$ satisfies
their conditions. We have not explored this question further, but we suspect
there is something to be found there.

\subsubsection*{AI use statement}

GPT-6 has been used to proofread the paper; it has also contributed Remark
\ref{rmk.not-enough} and drawn Figure~\ref{fig.sep.proof-dag}. All writing and
all proofs are our own.


\begin{thebibliography}{999999999}                                                                                        %


\bibitem[BjoBre05]{BjoBre}\href{https://doi.org/10.1007/3-540-27596-7}{Anders
Bjorner, Francesco Brenti, \textit{Combinatorics of Coxeter Groups}, Springer
2005}.\newline See \url{https://www.mat.uniroma2.it/~brenti/correct.ps} for errata.

\bibitem[Bourba02]{Bourbaki-Lie4-6}Nicolas Bourbaki, \textit{Lie Groups and
Lie Algebras, Chapters 4--6}, Springer 2002.

\bibitem[Davis25]{Davis}%
\href{https://people.math.osu.edu/davis.12/second_edition.pdf}{Michael W.
Davis, \textit{The Geometry and Topology of Coxeter Groups}, 2nd edition
2025.}\newline See
\url{https://www.cip.ifi.lmu.de/~grinberg/algebra/davis-gtcox-errata.pdf} for
unofficial errata.

\bibitem[Digne08]{Digne}Fran\c{c}ois Digne, \textit{Alg\`{e}bres de Hecke}, 5
December 2008.\newline\url{https://www.lamfa.u-picardie.fr/digne/hecke.pdf}

\bibitem[GaeGao20a]{GaeGao20a}%
\href{https://arxiv.org/abs/1905.09331v2}{Christian Gaetz, Yibo Gao,
\textit{Separable elements in Weyl groups}, arXiv:1905.09331v2, Adv. Appl.
Math. \textbf{113} (2020).}

\bibitem[GaeGao20b]{GaeGao20b}%
\href{https://arxiv.org/abs/1911.11172v2}{Christian Gaetz, Yibo Gao,
\textit{Separable elements and splittings of Weyl groups}, arXiv:1911.11172v2,
Adv. Math. \textbf{374} (2020).}

\bibitem[Gaetz21]{Gaetz21}%
\href{https://hdl.handle.net/1721.1/139101}{Christian Gaetz, \textit{New
combinatorics of the weak and strong Bruhat orders}, PhD thesis at MIT, 2021.}

\bibitem[Ghys16]{Ghys16}\href{https://arxiv.org/abs/1612.06373v4}{\'{E}tienne
Ghys, \textit{A Singular Mathematical Promenade}, arXiv:1612.06373v4}.\newline
Published in: ENS \'{e}ditions, Lyon, 2017.

\bibitem[Hiller82]{Hiller}Howard Hiller, \textit{Geometry of Coxeter groups},
Pitman Advanced Publishing Program 1982.

\bibitem[Humphr72]{Humphreys72}James E. Humphreys, \textit{Introduction to Lie
Algebras and Representation Theory}, Springer 1972.

\bibitem[Humphr90]{Humphreys}James E. Humphreys, \textit{Reflection groups and
Coxeter groups}, Cambridge University Press 1990.

\bibitem[Kitaev11]{Kitaev11}%
\href{https://doi.org/10.1007/978-3-642-17333-2}{Sergey Kitaev,
\textit{Patterns in Permutations and Words}, Springer 2011.}

\bibitem[Korotk16]{Korotk16}%
\href{https://doi.org/10.1007/s10958-017-3413-5}{S. Korotkikh, \textit{The
Mallows measures on the hyperoctahedral group}, Zapiski Nauchnykh Seminarov
POMI \textbf{448} (2016), pp. 151--164.}

\bibitem[Luszti14]{Lusztig}\href{http://arxiv.org/abs/math/0208154v2}{George
Lusztig, \textit{Hecke algebras with unequal parameters},
arXiv:math/0208154v2.}

\bibitem[Heckma18]{Heckman}Gert Heckman, \textit{Coxeter groups}, 10 October
2018.\newline\url{https://www.math.ru.nl/~heckman/CoxeterGroups.pdf}

\bibitem[Stanle24]{Stanley-EC2}Richard P. Stanley, \textit{Enumerative
Combinatorics, volume 2}, Second edition, Cambridge University Press
2024.\newline See \url{http://math.mit.edu/~rstan/ec/} for errata.

\bibitem[West95]{West95}%
\href{https://doi.org/10.1016/0012-365X(94)00067-1}{Julian West,
\textit{Generating trees and the Catalan and Schr\"{o}der numbers}, Discrete
Math. \textbf{146} (1995), pp. 247--262.}
\end{thebibliography}
\end{document}